\documentclass[12pt,draftcls,onecolumn]{IEEEtran}
\usepackage{dsfont}
\usepackage{amssymb}
\usepackage{amsmath}
\usepackage[dvips]{graphicx}
\newtheorem{thm}{Theorem}
\newtheorem{lem}{Lemma}
\newtheorem{rmk}{Remark}
\newtheorem{defn}{Definition}
\newtheorem{prop}{Proposition}
\newtheorem{cor}{Corollary}
\newtheorem{ex}{Example}
\newtheorem{problem}{Problem}

\def \H {\mathcal{H}}
\def\<{\langle}
\def\>{\rangle}

\begin{document}
%
\title{Ensemble Control of Finite Dimensional Time-Varying Linear Systems}
%
%
\author{Jr-Shin~Li
\thanks{
        This work was supported by the NSF SGER and CAREER grants.}
\thanks{J.-S. Li is with the Department of Electrical and Systems Engineering, Washington University, St. Louis, MO 63130 USA (e-mail: jsli@seas.wustl.edu).}}

\maketitle

\begin{abstract}
In this article, we investigate the problem of simultaneously steering an uncountable family of finite dimensional time-varying linear systems. We call this class of control problems Ensemble Control, a notion coming from the study of spin dynamics in Nuclear Magnetic Resonance (NMR) spectroscopy and imaging (MRI). This subject involves controlling a continuum of parameterized dynamical systems with the same open-loop control input. From a viewpoint of mathematical control theory, this class of problems is challenging because it requires steering a continuum of dynamical systems between points of interest in an infinite dimensional state space by use of the same control function. The existence of such a control raises fundamental questions of ensemble controllability. We derive the necessary and sufficient controllability conditions and an accompanying analytical optimal control law for ensemble control of time-varying linear systems. We show that ensemble controllability is in connection with singular values of the operator characterizing the system dynamics. In addition, we study the problem of optimal ensemble control of harmonic oscillators to demonstrate our main results. We show that the optimal solutions are pertinent to the study of time-frequency limited signals and prolate spheroidal wave functions. A systematic study of ensemble control systems has immediate applications to systems with parameter uncertainties as well as to broad areas of quantum control systems as arising in coherent spectroscopy and quantum information processing. The new mathematical structures appearing in such problems are an excellent motivation for new developments in control theory.

\end{abstract}


\IEEEpeerreviewmaketitle

\section{Introduction}
\label{intro} State-of-the-art quantum technology can trap and experiment with individual atoms, image brains as well as generate structural and dynamical information of biological macromolecules. Numerous applications arising from such emerging techniques involve controlling a large quantum ensemble, e.g., on the order of Avogadro number ($6\times 10^{23}$), by use of the same control field \cite{Li_PRA, Li_preprint, Li_thesis, Li_CDC2006, Li_NOLCOS2007, Li_IEEE, Brent1, Brent2, Skinner1, Kobzar, Pattern}. In many cases, the elements of the ensemble show variations in the values of parameters characterizing the system dynamics, and hence the system Hamiltonians of different elements of the ensemble are distinct. This phenomenon results in a dispersion in the system dynamics. For example, in magnetic resonance experiments, the spins of an ensemble may have large dispersions in their natural frequencies (Larmor dispersion), strength of the applied radio frequency (RF) field (RF inhomogeneity), and the dissipation rates of the spins \cite{Levitt, Tycko, Clare}. In solid state NMR spectroscopy of powder, the random distribution of the orientations of internuclear vectors of coupled spins within an ensemble leads to a distribution of coupling strengths \cite{Rohr}.

A canonical problem among these applications is to develop excitations (control signals) that will steer such an ensemble of systems with different dynamics from an initial state to a desired final state in finite time using a time-varying electromagnetic pulse. From the perspective of mathematical control theory, this is a very challenging state transfer problem because it requires steering a continuum of dynamical systems between points of interest in an infinite dimensional state space with the same control function. This motivates the study of \textbf{Ensemble Control} and the notion of ensemble controllability \cite{Li_PRA, Li_preprint, Li_thesis, Li_CDC2006, Li_NOLCOS2007, Li_IEEE} described as follows.

Consider a parameterized family of control systems
\begin{align}
\label{eq:dxs}& \frac{d}{dt}X(t,s)=F(X(t,s), u(t), t, s),\\
& X\in M\subset\mathbb{R}^{\rm n},\quad s\in D\subset\mathbb{R}^{\rm d},\quad u\in U\subset\mathbb{R}^{\rm m},\nonumber
\end{align}
where $F$ is a smooth function of its arguments and $D$ is a compact subset of $\mathbb{R}^{\rm d}$. Different values of the parameter $s$ in (\ref{eq:dxs}) correspond distinct members of the ensemble $X(t,s)$ showing variations, but we are constrained to use the same open-loop control $u(t)$ to steer the whole ensemble. The existence of such a control raises fundamental questions of \emph{ensemble controllability}. The formal definition will be given in Section \ref{sec:basics}. In practice, such control designs are called \emph{compensating pulse sequences} as they can compensate for or are insensitive to the dispersion in system dynamics. Typical applications include the design of excitation and inversion pulses in NMR spectroscopy in the presence of Larmor dispersion and RF inhomogeneity \cite{Skinner1, Kobzar, Pattern, Levitt, Tycko, Levitt1, Garwood}, the transfer of coherence between a coupled spin ensemble with variations in the coupling strengths \cite{Chingas}, and the construction of slice selective pulses in MRI, where some spins of the ensemble are excited or inverted while the others remain unaffected \cite{Silver, Rourke, Shinnar, Roux, Conolly, Mao, Rosenfeld}. Practical considerations, such as power constraints and signal or information losses due to relaxation effects, make it desirable to construct pulses that achieve a desired level of compensation with minimum energy or in the shortest possible time. These considerations give rise to problems in \emph{optimal control of ensembles}. These pulse design problems are widely studied in NMR spectroscopy on the subject of composite pulses that correct the dispersion in system dynamics \cite{Levitt, Levitt1, Garwood, Tycko1, Shaka2, Levitt2}. However, a systematic study of the design of compensating pulse sequences has been missing. The research in Ensemble Control will give explicit answers to the questions of when a compensation for the system dynamics is possible and how, if possible, to achieve the desired level of compensation.

More generally, ensemble control provides a framework to devise open-loop controls that are robust to parameter uncertainties, namely, they are insensitive to parameters. In many control applications, an accurate model is not available and systems have unknown parameters.  They either have some level of uncertainty or are not deterministic, however, there are instead bounds or distributions that describe these parameters. For example, systems biology models have numerous parameters, such as kinetic constants, which are unknown or only weakly constrained by existing experimental knowledge \cite{Gutenkunst}. In chaotic dynamics, the problem of parameter uncertainties is unavoidable in synchronizing chaotic systems \cite{Etemadi}. In such scenarios, one aims to design controllers that are robust to these parameter uncertainties. Subjects on robust control theory and sliding mode control are well studied to design controllers that can control or stabilize such systems using ``feedback'' \cite{Byrnes, Edwards}. While effective, these controllers are closed-loop and dependent on measurement of the system state - at times a difficult, or impossible requirement. Ensemble control, rather, provides a systematic framework for the design of open-loop controls that are immune to parameters and robust to uncertainties.

This paper is organized as follows. In the following section, we review ensemble control problems and the notion of ensemble controllability. We summarize our previous work of steering an ensemble of systems evolving on SO(3) to highlight the basics of ensemble control. Simple examples of ensemble control of linear systems are addressed to motivate the need of developing ensemble controllability conditions. In Section \ref{sec:theory}, we present our main results of necessary and sufficient controllability conditions for ensemble control of finite-dimensional time-varying linear systems. Finally, we study optimal control of an ensemble of harmonic oscillators as a demonstration of our main results as well as provide an insight of how ensemble control framework can be adopted to deal with systems with parameter uncertainties. We show that this system is ensemble controllable and the optimal control of such a system is pertinent to the study of time-frequency limited signals and prolate spheroidal wave functions (pswf). We present both analytical and numerical solutions.

\section{Basics of Ensemble Control}
\label{sec:basics} In this section, we review the basics of the ensemble control and the notion of ensemble controllability \cite{Li_thesis, Li_IEEE}. Ensemble Control involves problems of simultaneously manipulating a continuum of dynamical systems with different internal and external dynamics by use of the same open-loop control input. The general form of an ensemble control system is shown as in (\ref{eq:dxs}), and standard linear and bilinear ensemble control systems are of the respective forms 
\begin{align*}
\frac{d}{dt}X(t,s) &= A(t,s)X(t,s)+B(t,s)u(t),\\
\frac{d}{dt}X(t,s) &= \Big[A(t,s)+\sum_{i}u_i(t) B_i(t,s)\Big]X(t,s).
\end{align*}
The investigation of what kind of dispersions $s$ in the system dynamics can and cannot be corrected is a subject of fundamental and practical importance. This raises interesting questions of ensemble controllability.

\begin{figure}[t]
\centering
\begin{tabular}{cc}
{\small (a)}\includegraphics[scale=0.7]{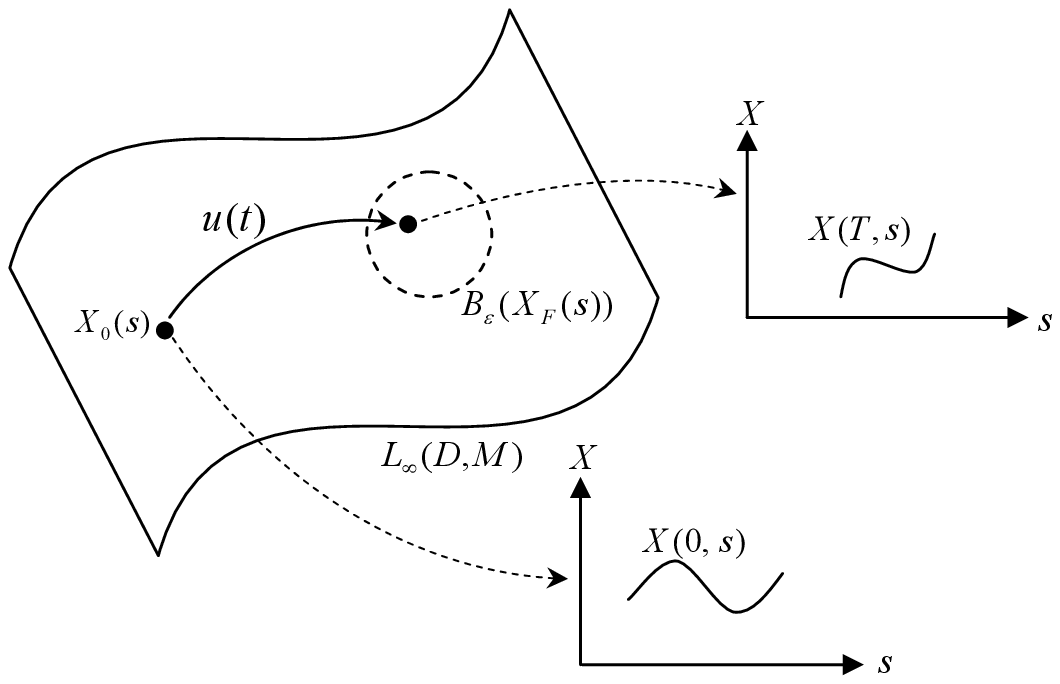}&\,
{\small (b)}\includegraphics[scale=0.70]{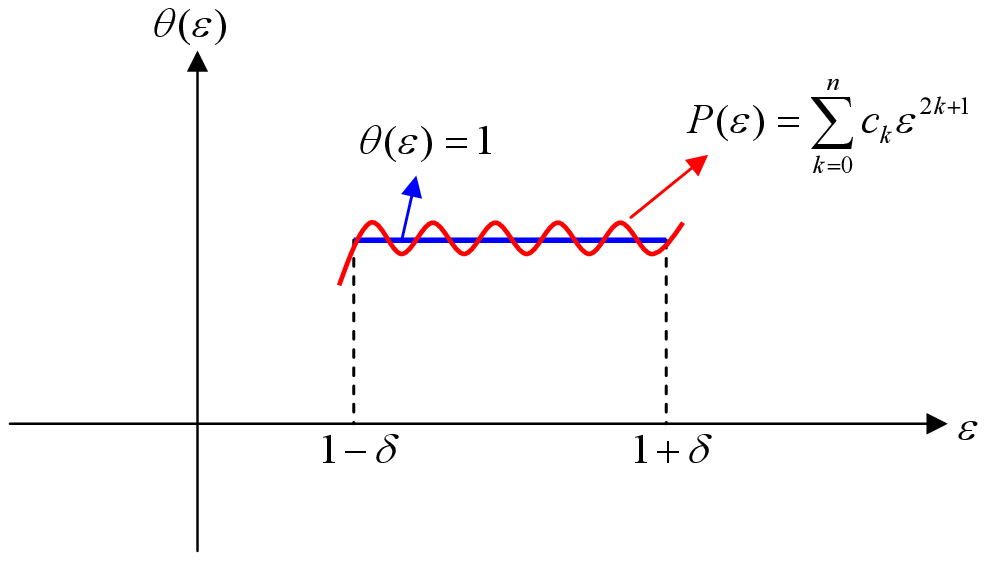}
\end{tabular}
\caption{\small (a) Illustration of the idea of ensemble controllability. $X_0(s)$ and $X(T,s)$ are two points on the function space $L_{\infty}(D,M)$, which correspond to two functions on the $s-X$ domain. If there exists a $u(t)$ that steers the system (\ref{eq:dxs}) from an initial point $X_0(s)$ to $X(T,s)\in\mathbf{B}_{\varepsilon}(X_F(s))$, for some finite time $T$, then the system is ensemble controllable. (b) The idea of the polynomial approximation, where the constant function $\theta(\epsilon)=1$ is approximated by an odd polynomial of degree $2n+1$, i.e., $\sum_{k=0}^{n}c_k\varepsilon^{2k+1}\approx 1$ for all $\varepsilon\in[1-\delta,1+\delta]$.}
\label{fig:ec}
\end{figure}

\begin{defn}
\label{def:controllability}
Consider a family of control systems as in (\ref{eq:dxs}). This family is called {\it ensemble controllable} on the function space $L_{\infty}(D,M)$ if and only if for all $\varepsilon>0$ and for all $X_0,X_F\in L_{\infty}(D,M)$, there exists $T>0$ and an open-loop piecewise-continuous control $u:[0,T]\rightarrow U$ such that starting from any initial state $X_0(s)=X(0,s)$, the final state $X_T(s)=X(T,s)\in L_{\infty}(D,M)$ satisfies $\|X_T-X_F\|_{\infty}\leq\varepsilon$.
\end{defn}

Note that $T\in(0,\infty)$ may depend on $\varepsilon$, $D$, and the bound of the control amplitude. The idea of ensemble controllability is illustrated in Figure \ref{fig:ec}(a).

\begin{rmk}
Note that ensemble control systems are different from distributed parameter systems where the systems are governed by partial differential equations \cite{Stavroulakis}. The difference can be seen from the fact that in (\ref{eq:dxs}) there is no partial derivative term with respect to the second variable, $\frac{\partial x}{\partial s}$. Ensemble control is also different from robust control in the sense that in ensemble control one is interested in devising an open-loop control signal, $u(t)$, that is insensitive to the state variable $x$ and the parameter $s$, while in most cases of robust control designs, constructing a closed-loop feedback control is of interest \cite{Feintuch}.
\end{rmk}

\subsection{The Prototype of Ensemble Control Problem}
The prototype of ensemble control problem was presented in our previous work, that is, control of a continuum of systems on SO(3) \cite{Li_thesis, Li_IEEE},
$$\dot{X}(t,\omega,\epsilon)=\Big[\omega\Omega_z+\epsilon u\Omega_y+\epsilon v\Omega_x\Big]X(t,\omega,\epsilon),\quad X(0,\omega,\epsilon)=I,$$ where $X\in SO(3)$, $(\omega,\epsilon)\in D=[a,b]\times[c,d]\in\mathbb{R}\times\mathbb{R}^{+}$, and $\Omega_x,\Omega_y,\Omega_z$ are generators of rotation around $x,y$, and $z$ axis, respectively. This system is ensemble controllable on the function space $\mathcal{S}(D)$, the set of all SO(3) valued measurable functions defined on $D$. We also showed through the study of this system that generating higher order Lie brackets by use of the control vector fields which carry higher order powers of the dispersion parameters is a key to investigating ensemble controllability. As a result, the infinite dimensional control problem can be translated to the problem of polynomial approximation. In this example, we can synthesize new generators of rotations by successive Lie bracketing of $\{\omega\Omega_z,\epsilon\Omega_y,\epsilon\Omega_x\}$ and then produce $(\omega,\epsilon)$-dependent evolutions of the form,
$$\exp\bigg\{\sum_{l}\Big(\sum_{k}c_{kl}\omega^{k}\Big)\epsilon^{2l+1}\Omega_x\bigg\}.$$
This evolution can then be used to approximate the desired evolution, $\exp\{\theta(\omega,\epsilon)\Omega_x\}$, with dependence on $(\omega,\epsilon)$ by appropriate choice of $c_{kl}$'s, where $\theta(\omega,\epsilon)$ is a measurable function over $D$.

The idea of polynomial approximation is shown in Figure \ref{fig:ec}(b). Based on this concept, a result on ensemble controllability of a simple linear system is immediately evident.
%
\begin{ex}
\rm An ensemble of time-invariant linear systems $$\dot{X}(t,s)=A X(t,s) + sBu(t),$$
where $s\in[s_1,s_2]\subset\mathbb{R}$, $A\in\mathbb{R}^{n\times n}$ and $B\in\mathbb{R}^{n\times m}$, is not ensemble controllable \cite{Li_IEEE}.\\
\end{ex}
{\it Proof.} Observe that $$sBu=\sum_{k=1}^m su_k b^k,$$ where $b^k$ is the $k$th column of $B$. We can think of $b^k$ as constant vector fields that generate translations. Since $b^k$ all commute to one another, their Lie brackets do not generate terms carrying higher powers of the dispersion parameter $s$. Therefore, the system is not ensemble controllable. Note that this result can also be easily verified by applying the variation of constants formula.\hfill $\Box$

The above example shows that the inability to synthesize higher powers of the dispersion parameter by successive Lie bracketing makes this system not ensemble controllable. The following fact characterizes a necessary controllability condition for a family of single-input linear systems.

\begin{thm}
Consider a family $(A(s),b(s))$ of linear single-input controllable systems
$$\dot{X}(t,s)=A(s)X(t,s)+b(s)u(t),$$
where $s$ takes values from a finite set $D\subset\mathbb{R}$, $A(s)\in\mathbb{R}^{n\times n}$, $b(s)\in\mathbb{R}^{n}$, and $u:[0,T]\rightarrow\mathbb{R}$, $T<\infty$. If the system is controllable, then there are no repeated eigenvalues of $A(s)$ for all $s$ \cite{Li_thesis}.
\end{thm}

\begin{rmk}
For the linear single-input controllable system $(A(s),b(s))$, it is required that $rank\ A(s)\geq n-1$, and the equality holds when $A(s)$ contains at least one eigenvalue equal to $0$ since $\det A(s)=0$. However, the condition that there are no repeated eigenvalues among all $A(s)$ restricts the ensemble to contain at most one singular $A(s)$.
\end{rmk}

These examples motivate the need of developing general controllability conditions. In the following section, we investigate necessary and sufficient controllability conditions for an ensemble of finite dimensional time-varying linear systems. We show that these conditions are associated with the singular value representation of the linear operator characterizing system dynamics.

\section{Linear Operators and Ensemble Controllability}
\label{sec:theory} The main result of this article is to provide the necessary and sufficient controllability conditions for an ensemble of general finite dimensional time-varying linear systems of the form
\begin{eqnarray}
\label{eq:linearTV}
\frac{d}{dt}X(t,s)=A(t,s)X(t,s)+B(t,s)u(t),
\end{eqnarray}
where $A(t,s)$ and $B(t,s)$ are $n\times n$ and $n\times m$ matrices, respectively, whose elements are complex-valued $L_2$ functions defined on a compact set $D=[0,T]\times [s_1,s_2]\subset\mathbb{R}^2$, denoted as $A\in L_2^{\rm n\times n}(D)$ and $B\in L_2^{\rm n\times m}(D)$. In this case, the ensemble controllability conditions are associated with when there exists an open-loop control, $u\in L_2^{\rm m}[0,T]$, which will steer the whole ensemble $X$ between points of interest in the function space $L_2^{\rm n}[s_1,s_2]$.

Let's start with some standard control theoretic analysis for the above system. Consider a fixed finite time $T$, starting from an initial state $X(0,s)$ we have by the variation of constants formula
\begin{eqnarray}
\label{eq:XT}
X(T,s)=\Phi(T,0;s)X(0,s)+\int_0^T \Phi(T,\tau; s)B(\tau,s)u(\tau)d\tau,
\end{eqnarray}
where $\Phi(t,0;s)$ is the transition matrix for $\frac{d}{dt}X(t,s)=A(t,s)X(t,s)$. It is known that for each $s\in[s_1,s_2]$
\begin{align*}
\Phi(t,0;s)=I &+ \int_{0}^{t}A(\sigma_1,s)d\sigma_1+\int_{0}^{t}A(\sigma_1,s) \int_{0}^{\sigma_1}A(\sigma_2,s)d\sigma_2 d\sigma_1\\
&+ \int_{0}^{t}A(\sigma_1,s)\int_{t_0}^{\sigma_1}A(\sigma_2,s)\int_{t_0}^{\sigma_2} A(\sigma_3,s) d\sigma_3 d\sigma_2 d\sigma_1+\ldots,
\end{align*}
is the Peano-Baker series which is uniformly convergent \cite{Brockett}. Given a desired target state $X_F(s)$ and an $\varepsilon>0$, we wish to find a control $u(t)$ such that $\|X(T,s)-X_F(s)\|_2\leq\varepsilon$. A simple manipulation of \eqref{eq:XT} yields
\begin{eqnarray}
\label{eq:xi}
\int_0^T\Phi(0,\tau;s)B(\tau,s)u(\tau)d\tau=\xi(s),
\end{eqnarray}
where $$\xi(s)=\Phi(0,T;s)X_F(s)-X(0,s),$$
and $\xi(s)$ is known as long as the initial and target states are specified. Note that for now we consider $X(T,s)=X_F(s)$. Let $\mathcal{H}_1=L_2^{\rm{m}}[0,T]$ be the set of $m$-tuples, whose elements are complex vector-valued square-integrable measurable functions defined on $0\leq t\leq T$, with an inner product defined by
\begin{eqnarray}
\label{eq:inner1}
\< g,h\> _{\H_1} = \int_0^T g^{\dagger}(t) h(t)dt,
\end{eqnarray}
where $\dagger$ denotes the conjugate transpose. Let $\H_2=L_2^{\rm{n}}[s_1,s_2]$ equipped with an inner product
\begin{eqnarray}
\label{eq:inner2}
\< p,q\> _{\H_2} = \int_{s_1}^{s_2} p^{\dagger}(s) q(s)ds.
\end{eqnarray}
It is clear that, with well-defined addition and scalar multiplication, $\H_1$ and $\H_2$ are separable Hilbert spaces. Now we define $L:\H_1\rightarrow\H_2$ by
$$(Lu)(s)=\int_0^T \Phi(0,\tau;s)B(\tau,s)u(\tau)d\tau,$$
and hence from (\ref{eq:xi})
\begin{eqnarray}
\label{eq:Lu}
(Lu)(s)=\xi(s).
\end{eqnarray}
Denote $\mathcal{B}(\H_1,\H_2)$ as the set of bounded linear operators from  $\H_1$ to $\H_2$. It is then easy to verify that $L\in\mathcal{B}(\H_1,\H_2)$ (see Appendix \ref{apd:Lcompact}). Consequently, $L$ has the adjoint $L^*$ satisfying
$$\left \< f,Lu\right \> _{\H_2}=\left \< L^*f,u \right \> _{\H_1},\quad\forall f\in\H_2, u\in\H_1.$$
This gives, by (\ref{eq:inner1}), (\ref{eq:inner2}), and the Fubini's theorem,
\begin{align*}
\int_{s_1}^{s_2}f^{\dagger}(s)\left(\int_0^{T}\Phi(0,\tau;s)B(\tau,s)u(\tau)d\tau\right)ds &=\int_0^T\left(\int_{s_1}^{s_2}\left[B^{\dagger}(\tau,s)\Phi^{\dagger}(0,\tau;s)f(s)\right]^{\dagger}ds\right)u(\tau)d\tau\\
&=\int_0^T\left(L^* f\right)^{\dagger}u(\tau)d\tau.
\end{align*}
Therefore,
\begin{eqnarray}
\label{eq:L*}
(L^*f)(t)=\int_{s_1}^{s_2}B^{\dagger}(t,s)\Phi^{\dagger}(0,t;s)f(s)ds.
\end{eqnarray}
The study of ensemble controllability for this system boils down to the problem of solving the inverse problem as in (\ref{eq:Lu}). With the above analysis, we show the main result of this paper.

\begin{thm}
\label{thm:main}
Consider a parameterized family of finite dimensional time-varying linear systems
\begin{eqnarray}
\label{eq:maineq}
\frac{d}{dt}X(t,s)=A(t,s)X(t,s)+B(t,s)u(t),
\end{eqnarray}
where $X:D\rightarrow\mathbb{R}^{n}$, $D=[0,T]\times[s_1,s_2]\subset\mathbb{R}^2$, and $u\in L_2^{\rm m}$; $A\in L_2^{\rm n\times n}(D)$, $B\in L_2^{\rm n\times m}(D)$, and time-varying $(A,B)$ are controllable pairs for all $s$. This family is ensemble controllable on the function space $L_2^{\rm n}[s_1,s_2]$ if and only if
\begin{align}
{\rm (i)} & \quad\sum_{n=1}^{\infty}\frac{|\< \xi,\nu_n\>|^2}{\sigma_n^2}<\infty\\
{\rm (ii)} & \quad\xi\in\overline{\mathcal{R}(L)},
\end{align}
where $(\sigma_n,\mu_n,\nu_n)$ is a singular system of $L$. Moreover, the control law
$$u=\sum_{n=1}^{\infty}\frac{1}{\sigma_n}\< \xi,\nu_n\> \mu_n,$$
satisfies $$\<u,u\>\leq\<u_0,u_0\>$$ for all $u_0\in\mathcal{U}$ and $u_0\neq u$, where $\mathcal{U}=\big\{v\,|\,Lv=\xi\ \text{with (i) and (ii)}\big\}$. In addition, $$u_N=\sum_{j=1}^{N(\varepsilon)}\frac{1}{\sigma_j}\< \xi,\nu_j\> \mu_j,$$
is the best approximation of $u$ for a given $\varepsilon>0$, namely, $u_N$ is such that $\|\xi-Lu_m\|\leq\varepsilon$ for all $m\geq N(\varepsilon)$, where $$u_m=\sum_{j=1}^{m}\frac{1}{\sigma_j}\< \xi,\nu_j\> \mu_j.$$ \hfill$\blacksquare$
\end{thm}

Before proving the above theorem, we need the following preliminary tools.

\subsection{Preliminaries}
\begin{prop}
\label{prop:Lcompact}
The operator $L:\H_1\rightarrow\H_2$ defined by
$$(Lu)(s)=\int_0^T \Phi(0,\tau;s)B(\tau,s)u(\tau)d\tau,$$ is compact.
\end{prop}
{\it Proof.} See Appendix \ref{apd:Lcompact}.

\begin{thm}[Spectral Theorem \cite{Porter}]\label{thm:spectral}
Let $X$ be a Hilbert space and $A:X\rightarrow X$ be a compact self adjoint operator. Then there exist a, possibly finite, sequence $\{\mu_n\}$ of nonzero eigenvalues of $K$ and a corresponding orthonormal sequence $\{\phi_n\}$ of eigenvectors such that for each $x\in X$, $Ax=\sum_n\mu_n\<x,\phi_n\>\phi_n$, where the sum is a finite sum if there are only finitely many eigenvalues. Moreover if $\{\mu_n\}$ is an infinite sequence, then it converges to zero.
\end{thm}

\begin{defn}[Singular System \cite{Gohberg}] 
\label{def:singular}
Let $X$ and $Y$ be Hilbert spaces and $K:X\rightarrow Y$ be a compact operator. If $(\sigma_n^2,\nu_n)$ is an eigensystem of $KK^*$ and $(\sigma_n^2, \mu_n)$ is an eigensystem of $K^*K$, namely, $KK^* \nu_n=\sigma_n^2 \nu_n$ and $K^* K \mu_n=\sigma_n^2\mu_n$, where $\sigma_n>0$ ($n\geq 1$), and the two systems are related by the equations
\begin{eqnarray}
\label{eq:singular}
K \mu_n=\sigma_n \nu_n,\quad K^* \nu_n=\sigma_n\mu_n,
\end{eqnarray}
we say that $(\sigma_n,\mu_n,\nu_n)$ is a singular system of $K$.
\end{defn}

\begin{rmk}
Since $K$ is compact and then we know $KK^*$ and $K^*K$ are both compact, self-adjoint, and nonnegative operators. Thus by the Spectral theorem, $K^*K$ can be represented in terms of its positive eigenvalues, namely, we have $K^*Kx=\sum_{n}\sigma_n^2\<x,\mu_n\>\mu_n$ for all $x\in X$. Moreover, since $K^*K\mu_n=\sigma_n^2\mu_n$, the relations as in \eqref{eq:singular} follow by taking $\nu_n=\frac{1}{\sigma_n}K\mu_n$. This can be treated as the infinite dimensional analogue of the singular value decomposition of a matrix.
\end{rmk}

The above definition immediately gives rise to the following results.

\begin{prop}
\label{prop:singular}
Let $X$ and $Y$ be Hilbert spaces and $K:X\rightarrow Y$ be a compact operator. If $(\sigma_n,\mu_n,\nu_n)$ is a singular system of $K$, then
\begin{enumerate}
\item[(i)] $\{\mu_n\}$ is an orthonormal basis of $\overline{\mathcal{R}(K^*)}$,
\item[(ii)] $\{\nu_n\}$ is an orthonormal basis of $\overline{\mathcal{R}(K)}$.
\end{enumerate}
\end{prop}
{\it Proof.} (i) Since $\mu_n=\frac{1}{\sigma_n^2}K^*K\mu_n\in\mathcal{R}(K^*K)$ and $K^*K$ is compact and self adjoint, the Spectral theorem, $K^*K x=\sum_j\sigma_j^2\<x,\mu_j\>\mu_j$, for all $x\in X$, implies that $span\{\mu_n\}$ is dense in $\mathcal{R}(K^*K)$. It follows that $\overline{span}\{\mu_n\}=\overline{\mathcal{R}(K^*K)}=\overline{\mathcal{R}(K^*)}$. (ii) can be proved similarly. \hfill$\Box$

\begin{thm}[Singular value expansion \cite{Gohberg}]\label{thm:singular}
Let $X$ and $Y$ be Hilbert spaces, $K:X\rightarrow Y$ be a compact operator and $\{(\sigma,\mu_n,\nu_n)\ |\ n\in\Delta\}$ be a singular system for $K$. Then
$$Kx=\sum_{n\in\Delta}\sigma_n\<x,\mu_n\>\nu_n,\quad K^* y=\sum_{n\in\Delta}\sigma_n\<y,\nu_n\>\mu_n,$$
for all $x\in X$, $y\in Y$. In particular, if
$$K_nx=\sum_{j=1}^{n}\sigma_j\<x,\mu_j\>\nu_j,\quad x\in X,$$
and $K$ is of infinite rank, namely, $\Delta=\mathbb{N}$, then
$$\|K-K_n\|\leq\sup_{j>n}\sigma_j\rightarrow 0 \quad as\quad n\rightarrow\infty.$$
\end{thm}
{\it Proof.} Since, by Proposition \ref{prop:singular}, $\{\nu_n\}$ is an orthonormal basis of $\overline{\mathcal{R}(K)}$, the Fourier expansion gives for all $x\in X$
$$Kx=\sum_{j}\<Kx,\nu_j\>\nu_j=\sum_j\<x,K^*\nu_j\>\nu_j=\sum_j\sigma_j\<x,\mu_j\>\nu_j.$$
The other part can be shown similarly. Also, we have
\begin{align*}
\|(K-K_n)\,x\|^2=\|\sum_{j>n}\sigma_j\<x,\mu_j\>\nu_j\|^2=\sum_{j>n}|\sigma_j|^2|\<x,\mu_j\>|^2
\leq\sup_{j>n}\ \sigma_j^2\ \|x\|^2.
\end{align*}
Therefore, $$\|K-K_n\|\leq\sup_{j>n}\ \sigma_j\rightarrow 0\quad as\quad n\rightarrow\infty,$$
since $\sigma_n^2$ is an eigenvalue of $K^* K$ and $\sigma_n^2\rightarrow 0$ as $n\rightarrow\infty$. \hfill$\Box$

\begin{thm}[Riesz-Fischer Theorem \cite{Gohberg}]\label{thm:Riesz}
Let $\{u_1,u_2,\ldots\}$ be an orthonormal set in a Hilbert space $X$ and let $\{\alpha_n\}$ be a sequence of scalars. Then
$$\sum_{n=1}^{\infty}|\alpha_n|^2 \quad\text{converges if and only if} \quad\sum_{n=1}^{\infty}\alpha_n u_n \quad\text{converges},$$
and, in that case,
$$\alpha_n=\< x,u_n\> \quad\forall\, n\in\mathbb{N}, \quad\text{where}\quad x=\sum_{n=1}^{\infty}\alpha_n u_n.$$
\end{thm}

\begin{thm}[Minimum Norm \cite{Luenberger}]
\label{thm:minnorm} Let $G$ and $H$ be Hilbert spaces and let $A\in\mathcal{B}(G, H)$ with range closed in $H$. Then, the vector $x$ of minimum norm satisfying $Ax=y$ is given by $x=A^{*}z$ where $z$ is any solution of $AA^{*}z=y$ and $A^{*}$ is the adjoint operator of $A$.
\end{thm}

{\it Proof of Theorem \ref{thm:main}.} {\it Necessity:} Since $L$ is compact by Proposition \ref{prop:Lcompact}, we then let $(\sigma_n,\mu_n,\nu_n)$ be a singular system of $L$. Now suppose that $u(t)$ is a solution to (\ref{eq:Lu}) and $u\in\H_1$. Then
$$\< \xi,\nu_n\> =\< Lu,\nu_n\> =\< u,L^* \nu_n\> =\sigma_n\< u,\mu_n\> .$$
Hence, $\< u,\mu_n\> =\frac{1}{\sigma_n}\< \xi,\nu_n\> .$ By the application of Bessel's inequality to $u$ and the orthonormal system $\{\mu_n\}$, we have
$$\sum_{n=1}^{\infty}\frac{|\< \xi,\nu_n\> |^2}{\sigma_n^2}=\sum_{n=1}^{\infty}|\< u,\mu_n\> |^2\leq\|u\|^2<\infty.$$
Furthermore, for any $\eta\in\mathcal{N}(L^*)$, i.e., $L^*\eta=0$ and $\eta\in\H_2$, we have
$$\< \xi,\eta\> =\< Lu,\eta\> =\< u,L^*\eta\> =0.$$
Therefore, $\xi\in\mathcal{N}(L^*)^{\perp}=\overline{\mathcal{R}(L)}.$

{\it Sufficiency:} Conversely, we suppose that both conditions (i) and (ii) are satisfied. Let $\alpha_n=\frac{1}{\sigma_n}{\<\xi,\nu_n\>}$, hence $\sum_{n=1}^{\infty}|\alpha_n|^2$ converges according to the condition (i). By the Riesz-Fischer theorem, there exists a $u\in\H_1$ so that
\begin{eqnarray}
\label{eq:u_RF}
u=\sum_{n=1}^{\infty}\alpha_n \mu_n,
\end{eqnarray}
and then
\begin{eqnarray}
\label{eq:alphan}
\alpha_n=\<u,\mu_n\>=\frac{1}{\sigma_n}\<\xi,\nu_n\>.
\end{eqnarray}
Note that $u\in\mathcal{N}(L)^{\perp}\subset\H_1$ since $\{\mu_n\}$ spans $\overline{\mathcal{R}(L^*)}$. Hence, from Theorem \ref{thm:singular} and (\ref{eq:alphan}), we obtain 
\begin{eqnarray}
\label{eq:Lu1}
Lu=\sum_{n=1}^{\infty}\sigma_n\<u,\mu_n\>\nu_n=\sum_{n=1}^{\infty}\<\xi,\nu_n\>\nu_n.
\end{eqnarray}
Since $\xi\in\overline{\mathcal{R}(L)}$ by the condition (ii) and $\{\nu_n\}$ spans $\overline{\mathcal{R}(L)}$ by Proposition \ref{prop:singular}, $\xi$ can be expressed by the Fourier expansion, $$\xi=\sum_{n=1}^{\infty}\<\xi,\nu_n\>\nu_n.$$
Combining this with \eqref{eq:Lu1} and \eqref{eq:u_RF}, we conclude that
\begin{eqnarray}
\label{eq:control}
u=\sum_{n=1}^{\infty}\frac{1}{\sigma_n}\<\xi,\nu_n\>\mu_n,
\end{eqnarray}
$u\in\mathcal{N}(L)^{\perp}\subset\H_1$, is a solution of (\ref{eq:Lu}).
We now put
$$u_N=\sum_{j=1}^{N}\frac{1}{\sigma_j}\<\xi,\nu_j\>\mu_j,$$
where $N\in\mathbb{N}$. By the fact that $\{\mu_n\}$ is an orthonormal sequence, we have
$$\|u-u_N\|^2=\sum_{j=N+1}^{\infty}\frac{1}{\sigma_j^2}\,|\<\xi,\nu_j\>|^2\rightarrow 0 \quad as \quad N\rightarrow\infty,$$
and then, by Theorem \ref{thm:singular},
$$\|Lu-Lu_N\|^2=\sum_{j=N+1}^{\infty}\sigma_j^2\,|\<u,\mu_j\>|^2\rightarrow 0 \quad as \quad N\rightarrow\infty.$$
Therefore, given any $\varepsilon>0$, we can find $u_N$ such that $\|\xi-Lu_N\|\leq\varepsilon$ for an appropriate choice of $N=N(\varepsilon)$. Moreover, since $u\in\mathcal{N}(L)^{\perp}=\overline{\mathcal{R}(L^*)}$, according to Theorem \ref{thm:minnorm}, $u$ is of minimum norm satisfying (\ref{eq:Lu}), that is,
$$\<u,u\>\leq\<u_0,u_0\>$$
for all $u_0\in\mathcal{U}$ and $u_0\neq u$, where $\mathcal{U}=\{v\,|\,Lv=\xi\ \text{with (i) and (ii)}\}$.  $\hfill\Box$


\begin{rmk}
The controllability condition (i) implies that an ensemble control law exists if and only if the Fourier coefficients $\<\xi,\nu_n\>$ with respect to the singular functions $\nu_n$ decay fast enough relative to the singular values $\sigma_n$. Note that $(\sigma_n^2,\nu_n)$ is an eigensystem of the Hermitian operator $LL^*:\H_2\rightarrow\H_2$ defined by
$$(LL^* z)(s)=\int_{s_1}^{s_2}\int_0^T\Phi(0,\tau;s)B(\tau,s)B^{\dagger}(\tau,\sigma)\Phi^{\dagger}(0,\tau;\sigma)z(\sigma) \,d\tau d\sigma,$$
an analogy of the classical controllability Gramian. The ensemble controllability condition (i) coincides with the so called Picard criterion in the literature of integral equations.

\end{rmk}


\begin{cor}
For any given initial state $X_0(s)=X(0,s)$, the ensemble control law $u$ as in (\ref{eq:control}) does not depend continuously on the target state $X_F(s)$.
\end{cor}
{\it Proof.} Suppose that $\delta_n$ is a perturbation of $\xi$ with $\delta_n\rightarrow 0$ as $n\rightarrow\infty$ and that $\tilde{\xi}_n=\xi+\delta_n$. Let $u$ and $\tilde{u}_n$ are solutions to the integral equations $Lu=\xi$ and $L\tilde{u}_n=\tilde{\xi}_n$, respectively. We now consider a qualified perturbation $\delta_n=a\sqrt{\sigma_n}\nu_n$, where $a\in\mathbb{R}$. It is clear that $\delta_n\rightarrow 0$ as $n\rightarrow\infty$ because $\{\nu_n\}$ is an orthonormal basis and $\sigma_n\rightarrow 0$ as $n\rightarrow\infty$. Thus, we have
$$\|\tilde{\xi}_n-\xi\|=|a|\sqrt{\sigma_n}\rightarrow 0 \quad {\it as} \quad n\rightarrow\infty.$$
However,
$$\tilde{u}_n=\sum_{n=1}^{\infty}\frac{1}{\sigma_n}\<\tilde{\xi}_n,\nu_n\>\mu_n=u+\frac{a}{\sqrt{\sigma_n}}\ \mu_n,$$
and hence
$$\|\tilde{u}_n-u\|=\frac{|a|}{\sqrt{\sigma_n}}\rightarrow\infty\quad {\it as} \quad n\rightarrow\infty.$$
Therefore, the control $u$ doesn't depend continuously on $\xi(s)$ and thus neither on the target state $X_F(s)=\Phi(T,0;s)\left[\xi(s)+X_0(s)\right]$. $\hfill\Box$


\section{Optimal Control of An Ensemble of Harmonic Oscillators}
\label{sec:optimal}
In this section, we study in detail the ensemble control of a family of harmonic oscillators that demonstrates our main results in Section \ref{sec:theory}. We show that this system is ensemble controllable and derive an analytical optimal control law. The analysis of this ensemble control system is related to the study of time-frequency limited signals and prolate spheroidal wave functions. Alternatively, this problem can be viewed as control of a harmonic oscillator with parameter uncertainty, where the frequency is unknown but only its range is provided.


\subsection{Unconstrained Optimal Ensemble Control}
We first look at a fixed end-point optimal ensemble control problem without constraints on the control signals.

\begin{problem}
\label{prob:harmonic1} \rm Consider an ensemble of harmonic oscillators with a variation in their natural frequencies
\begin{equation}
\label{eq:harmonic}
\frac{d}{dt}\left[\begin{array}{c}x(t,\omega)\\y(t,\omega)\end{array}\right] =\left[\begin{array}{cc}0&-\omega\\ \omega&0\end{array}\right] \left[\begin{array}{c}x(t,\omega)\\y(t,\omega)\end{array}\right]
+\left[\begin{array}{c}u(t) \\ v(t)\end{array}\right],
\end{equation}
where $\omega\in D=[\omega_1,\omega_2]\subset\mathbb{R}$, $X(\cdot,\omega)=(x(\cdot,\omega),y(\cdot,\omega))^T \in L_2^{\rm 2}(D)$, and $U=(u,v)^T \in L_2^{\rm 2}[0,T]$. Find controls $u(t)$ and $v(t)$ that steer this continuum of systems from an initial state $X_0=(x(0,\omega),y(0,\omega))^{T}$ to within a ball of radius $\varepsilon$ around the final state $X_F=(x_F(\omega), y_F(\omega))^{T}$ at time $T<\infty$, and minimize the cost functional
\begin{eqnarray}
\label{eq:J} J=\int_0^T \left[u(t)^2+v(t)^2\right]dt.
\end{eqnarray}
\end{problem}

We first observe that each element of the ensemble in (\ref{eq:harmonic}) with a frequency $\omega\in D$ is controllable, because the Gramian matrix is of full rank, i.e.,
$$rank \Big[B\Big|AB\Big]=2,$$
where $$A=\left[\begin{array}{cc}0&-\omega\\ \omega&0\end{array}\right],\quad B=\left[\begin{array}{cc}1 & 0 \\ 0 & 1\end{array}\right].$$
%
%
\begin{thm}
\label{thm:smallcontrol} An ensemble of harmonic oscillators modeled as in (\ref{eq:harmonic}) is ensemble controllable on $L_2^2(D)$. 
\end{thm}
{\it Proof.} Without loss of generality, we consider the frequency distributes in a symmetric domain $D_s=[-\beta,\beta]$ since the system (\ref{eq:harmonic}) with $\omega\in[\omega_1,\omega_2]$ can be readily transformed to a frame with $\omega\in D_s$ by a simple change of coordinate. Let's rewrite (\ref{eq:harmonic}) as $\dot{X}=\omega\Omega X+BU$, where
$$\Omega=\left[\begin{array}{cc}0&-1\\ 1&0\end{array}\right].$$
Let $\tilde{X}(t)=\exp(-\tilde{\omega}\Omega t)X$, then we obtain the differential equation in the new coordinate, $\dot{\tilde{X}}=(\omega-\tilde{\omega})\Omega\tilde{X}+\exp(-\tilde{\omega}\Omega t)BU$, with frequencies $\nu=\omega-\tilde{\omega}\in[\omega_1-\tilde{\omega},\omega_2-\tilde{\omega}]$. Taking $\tilde{\omega}=(\omega_1+\omega_2)/2$, then we have
$$\frac{d}{dt}\tilde{X}(t,\nu)=\left[\begin{array}{cc}0&-\nu\\ \nu&0\end{array}\right]\tilde{X}(t,\nu) +\left[\begin{array}{c}\tilde{u}(t) \\ \tilde{v}(t)\end{array}\right],$$
where $\tilde{u}(t)=u(t)\cos(\tilde{\omega}t)+v(t)\sin(\tilde{\omega}t)$, $\tilde{v}(t)=-u(t)\sin(\tilde{\omega}t)+v(t)\cos(\tilde{\omega}t)$, $\nu\in[-\beta,\beta]$, and $\beta=(\omega_2-\omega_1)/2$.

Let
\begin{align*}
p(t,\omega) &= x(t,\omega)+iy(t,\omega),\\
\alpha(t) &= u(t)+iv(t),
\end{align*}
where $i=\sqrt{-1}$. The system (\ref{eq:harmonic}) can then be written as
\begin{eqnarray*}
\label{eq:maineqn} \dot{p}(t,\omega)=i\omega p(t,\omega)+\alpha(t),
\end{eqnarray*}
with $p(0,\omega)=x(0,\omega)+iy(0,\omega)$. By the variation of constants formula, we have at time $T$
\begin{equation}
\label{eq:pomega} p(T,\omega)=e^{i\omega T}p(0,\omega) +\int_{0}^{T}e^{i\omega(T-\tau)}\alpha(\tau)d\tau,
\end{equation}
for all $\omega\in D_s$. This gives
\begin{equation}
\label{eq:varnat}
\int_{0}^{T}e^{-i\omega\tau}\alpha(\tau)d\tau=e^{-i\omega T}p(T,\omega)-p(0,\omega)\doteq\xi(\omega).
\end{equation}
Let $\H_1=L_{2}[0, T]$ and $\H_2=L_{2}[-\beta,\beta]$ be Hilbert spaces over $\mathbb{C}$. Defining the linear operator $L: \H_1\rightarrow \H_2$
by
\begin{equation}
\label{eq:L}
(L\alpha)(\omega)=\int_{0}^{T}e^{-i\omega\tau}\alpha(\tau)d\tau\doteq\int_0^T k(\omega,\tau)\alpha(\tau)d\tau,
\end{equation}
we then have from (\ref{eq:varnat}) and (\ref{eq:L}) that
\begin{equation}
\label{eq:Lalpha} (L\alpha)(\omega)=\xi(\omega).
\end{equation}
Observe that $L$ is bounded since for every $f\in\H_1$, by the Cauchy-Schwartz inequality,
\begin{eqnarray*}
\|Lf\|_{\H_2}\leq\|(\sqrt{T}\ \|f\|_{\H_1})\|_{\H_2}=\sqrt{2BT}\|f\|_{\H_1}.
\end{eqnarray*}
Therefore, $L$ has the adjoint $L^*$ defined by
\begin{eqnarray}
\label{eq:L*}
(L^* g)(t)=\int_{-\beta}^{\beta}k(\omega,t)^{\dagger}g(\omega)d\omega
=\int_{-\beta}^{\beta}e^{i\omega t}g(\omega)d\omega.
\end{eqnarray}
Moreover, since $L: \H_1\rightarrow \H_2$ and $k(\omega,\tau)\in L_2([-\beta,\beta]\times [0,T])$, $L$ is a Hilbert-Schmidt operator on $\H_1$ and hence $L$ is compact (the compactness can also be shown following the proof of Proposition \ref{prop:Lcompact}). According to Theorem \ref{thm:minnorm}, the function $\alpha$ of minimum norm satisfying (\ref{eq:Lalpha}) is given by
\begin{equation}
\label{eq:controllaw} \alpha(t)=L^{*}z(\omega),
\end{equation}
where $z$ satisfies
\begin{equation}
\label{eq:zomega} (Wz)(\omega)=\xi(\omega),
\end{equation}
and the operator $W:\H_2\rightarrow\H_2$ is defined by
\begin{align}
(Wz)(\omega) &= (LL^{*}z)(\omega)=\int_{-\beta}^{\beta}\frac{e^{i(\omega^{\prime}-\omega)T}-1}
{i(\omega^{\prime}-\omega)}z(\omega^{\prime})d\omega^{\prime}\nonumber\\
\label{eq:LLstar}&= \int_{-\beta}^{\beta}2\pi e^{i\frac{T}{2}(\omega^{\prime}-\omega)} \Big[\frac{\sin\big(\frac{T}{2}(\omega-\omega^{\prime})\big)}
{\pi(\omega-\omega^{\prime})}\Big]z(\omega^{\prime})d\omega^{\prime}
\end{align}
for $\omega, \omega^{\prime}\in [-\beta,\beta]$. A simple change of variables converts (\ref{eq:LLstar}) into
\begin{eqnarray}
\label{eq:W2} (Wz)(\Omega)=\int_{-1}^{1}2\pi e^{ic(\Omega^{\prime}-\Omega)} \left[\frac{\sin[c(\Omega-\Omega^{\prime})]}{\pi(\Omega-\Omega^{\prime})}\right]z(\Omega^{\prime})d\Omega^{\prime},
\end{eqnarray}
in which $c=\frac{\beta T}{2}$, $\Omega=\frac{\omega}{\beta}$, $\Omega^{\prime}=\frac{\omega^{\prime}}{\beta}$, and $\Omega,\Omega^{\prime}\in [-1, 1]$. Notice that the term inside the bracket in (\ref{eq:W2}) is the kernel of the following integral equation:
\begin{eqnarray}
\label{eq:pswf}\int_{-1}^{1}\frac{\sin[c(\Omega-\Omega^{\prime})]}{\pi(\Omega-\Omega^{\prime})} \,\psi_{n}(\Omega^{\prime},c)\,d\Omega^{\prime}=\kappa_{n}(c)\, \psi_{n}(\Omega,c),
\end{eqnarray}
where the $n^{th}$ eigenfunction $\psi_{n}(\Omega,c)$ is the well-known \emph{prolate spheroidal wave function} (\emph{pswf}), and $\kappa_{n}(c)$ is the associated eigenvalue \cite{Percival, Flammer, Slepian1, Slepian2, Slepian3}, where $\kappa_n>0$ and $\kappa_n\rightarrow 0$ as $n\rightarrow\infty$. Consequently, the $n^{th}$ eigenfunction and the corresponding eigenvalue for $W$ as in (\ref{eq:W2}) can be easily represented in terms of $\psi_{n}$ and $\kappa_{n}$ by $\phi_{n}=e^{-i\omega\frac{T}{2}}\psi_{n}$, and $\lambda_{n}=2\pi\kappa_{n}$. Note that $\psi_{n}$'s are orthogonal and complete on $L_{2}[-1,1]$ \cite{Percival}. Since $L$ is compact, $W=LL^*$ is compact. It can then be spectral decomposed by the orthonormal basis $\{\tilde{\phi}_n\}$ applied to (\ref{eq:zomega}) 
\begin{eqnarray}
\label{eq:Wz} Wz=\sum_{n=1}^{\infty}\lambda_{n}\<z,\tilde{\phi}_{n}\>\tilde{\phi}_{n}=\xi,\quad \tilde{\phi}_{n}=e^{-i\omega\frac{T}{2}}\frac{\psi_{n}}{\|\psi_{n}\|},
\end{eqnarray}
and this sequence is being uniformly convergent by the spectral theorem \cite{Gohberg}. It is also clear that $\{\tilde{\phi}_{n}\}$ is an orthonormal basis of $W$. The solution of (\ref{eq:Wz}) takes the form
\begin{eqnarray*}
z=\sum_{n=1}^{\infty}\frac{1}{\lambda_{n}}\<\xi,\tilde{\phi}_{n}\>\tilde{\phi}_{n}.
\end{eqnarray*}
Finally, we show that the above series $z(\omega)$ can be truncated to $z_{N}(\omega)$ so that $\|Wz-Wz_N\|\rightarrow 0$ as $N\rightarrow\infty$. Then, we obtain the best approximation of the minimum energy control law $\alpha_N=L^* z_N$ by (\ref{eq:controllaw}).

\begin{lem}
Given any $\varepsilon>0$, there exists a finite series $z_N$,
\begin{eqnarray}
\label{eq:zN}
z_N=\sum_{n=1}^{N}\frac{1}{\lambda_{n}}\<\xi,\tilde{\phi}_{n}\>\tilde{\phi}_{n},
\end{eqnarray}
such that
$$\|Wz-Wz_N\|\rightarrow 0 \quad {\it as} \quad N\rightarrow\infty,$$
where $N=N(\varepsilon)$ depends on the choice of $\varepsilon$.
\end{lem}
{\it Proof.} By the orthonormality of $\{\tilde{\phi}_{n}\}$, we get
$$Wz_N=\sum_{n=1}^{N}\<\xi,\tilde{\phi}_{n}\>\tilde{\phi}_{n}.$$
Let $a_{n}=\<\xi,\tilde{\phi}_{n}\>$, and then we have
\begin{eqnarray}
\label{eq:error1}
\|Wz-Wz_{N}\|^{2}=\sum_{N+1}^{\infty}|a_{n}|^{2}.
\end{eqnarray}
Since, by the Bessel's inequality, $$\sum_{n=1}^{\infty}|a_{n}|^{2}<\|\xi\|^2<\infty,$$
the error in (\ref{eq:error1}) can be made in response to the desired $\varepsilon$ by the selection of $N=N(\varepsilon)$. \hfill $\blacksquare$

It follows from (\ref{eq:L*}), (\ref{eq:controllaw}) and (\ref{eq:zN}) that
\begin{eqnarray}
\label{eq:alphaN}\alpha_{N}(t)=\int_{-\beta}^{\beta}e^{i\omega t}\sum_{n=1}^{N(\varepsilon)} \frac{1} {\lambda_{n}}\<\xi,\tilde{\phi}_{n}\>\tilde{\phi}_{n}\,d\omega,
\end{eqnarray}
which will steer the system (\ref{eq:harmonic}) from $p(0,\omega)$ to within the ball $\mathbf{B_\varepsilon}\big(p(T,\omega)\big)$ at time $T$, where $\mathbf{B_\varepsilon}\big(p(T,\omega)\big)=\{h\in\H_2\,:\,\|p(T,\omega)-h(\omega)\|\leq\varepsilon\}$. In addition, $\alpha_N$ is the best approximation, for the given $\varepsilon>0$, of the control law $\hat{\alpha}=u+iv$ that minimizes the cost functional $J$ as in (\ref{eq:J}). \hfill$\Box$

We now have ensemble controllability for the system (\ref{eq:harmonic}), however, this result fails when either $u(t)$ or $v(t)$ is not available.

\begin{cor}
\label{rmk:harmonic} An ensemble of systems as in (\ref{eq:harmonic}) is not ensemble controllable if either $u(t)\equiv 0$ or $v(t)\equiv 0$.
\end{cor}
{\it Proof.} Without loss of generality, we suppose now that $v(t)=0$ and that the initial state $(x(0,\omega),y(0,\omega))=(0,0)$ for all $\omega\in D_s$. Note that each element of the ensemble is still controllable in this case. Let
\begin{eqnarray}
\tilde{X}(t,\omega)&=&x(t,\omega)-x(t,-\omega),\nonumber\\
\tilde{Y}(t,\omega)&=&y(t,\omega)+y(t,-\omega).\nonumber
\end{eqnarray}
The system described in (\ref{eq:harmonic}) can then be transformed to
\begin{equation}
\label{eq:harmonic1}
\frac{d}{dt}\left[\begin{array}{c}\tilde{X}\\ \tilde{Y}\end{array}\right] =\left[\begin{array}{cc}0&-\omega\\ \omega&0\end{array}\right]\left[\begin{array}{c}\tilde{X}\\ \tilde{Y}\end{array}\right]; \,\,\left[\begin{array}{c}\tilde{X}(0,\omega)\\ \tilde{Y}(0,\omega)\end{array}\right]
=\left[\begin{array}{c}0\\0\end{array}\right].\nonumber
\end{equation}
Since the above system is autonomous, it stays at the origin for all $t$, i.e., $(\tilde{X}(t,\omega),\tilde{Y}(t,\omega))\equiv (0,0)$. Thus, the system is not ensemble controllable. \hfill $\Box$

\begin{figure}[t]
\centering
\begin{tabular}{cc}
{\small (a)}\includegraphics[scale=0.53]{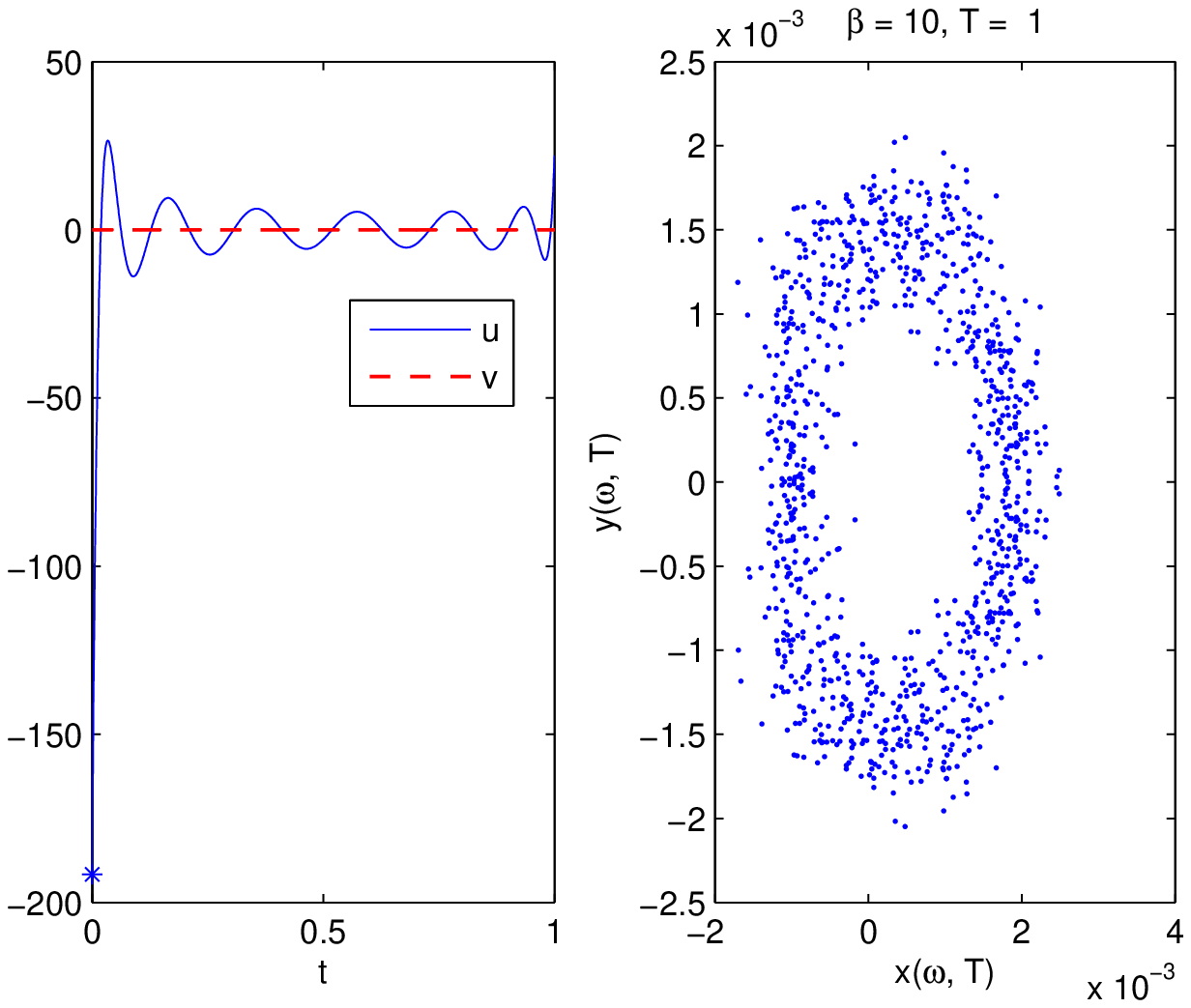}&
{\small (b)}\includegraphics[scale=0.53]{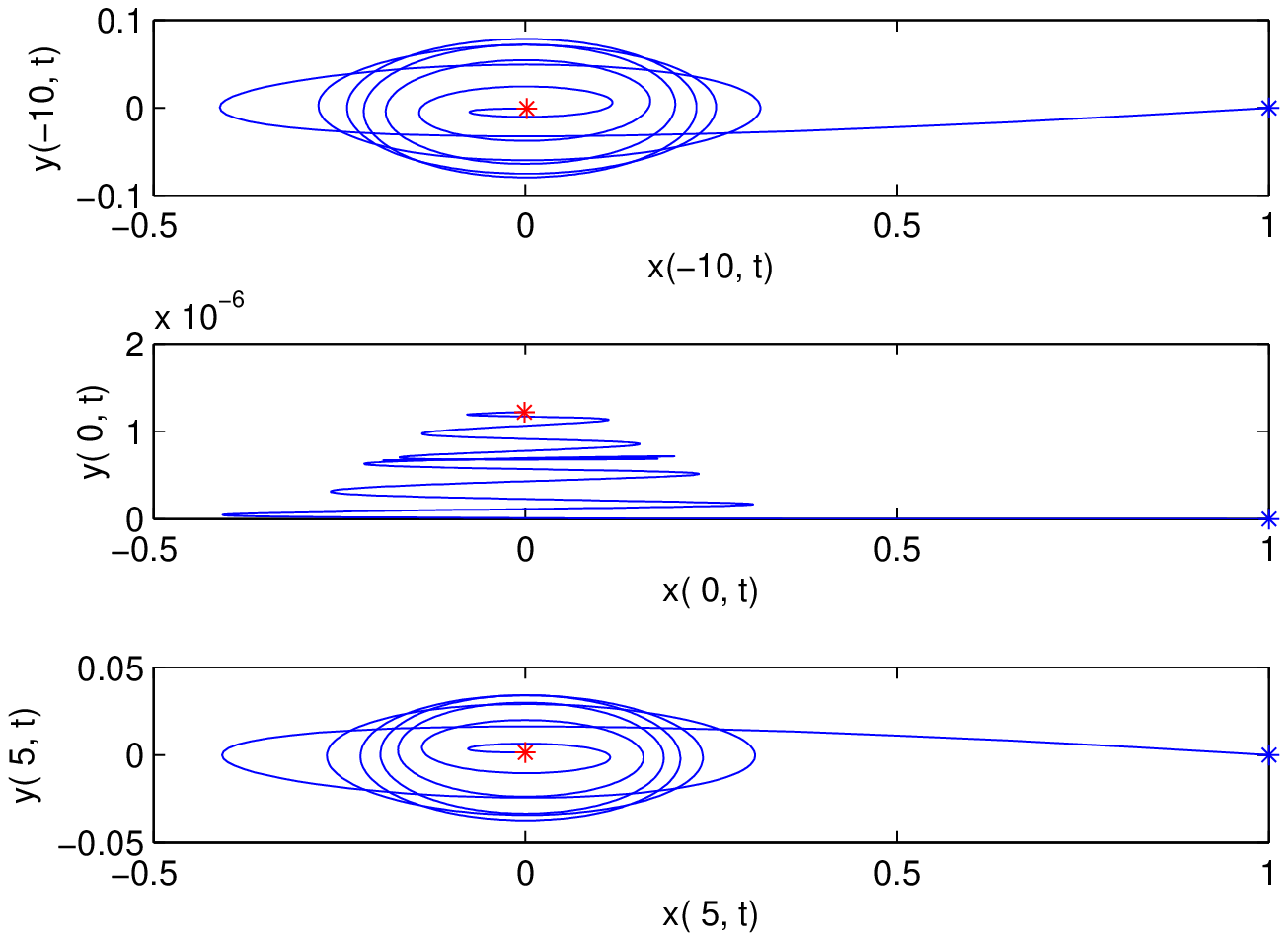}
\end{tabular}
\caption{\small The simulation results of Problem \ref{prob:harmonic1} for $N=1001$, $T=1$, and $\beta=10$. The initial state $X_0=(1,0)$ and the target state $X_F=(0,0)$. (a) The optimal control law $(u(t),v(t))$ for $t\in[0,1]$, and the final states for all systems $\omega\in[-10,10]$. (b) The trajectories for $\omega=-10$, $\omega=0$ and $\omega=5$ following $(u(t),v(t))$.}\label{fig:case1}
\end{figure}

\begin{figure}[h]
\centering
\begin{tabular}{cc}
{\small (a)}\includegraphics[scale=0.54]{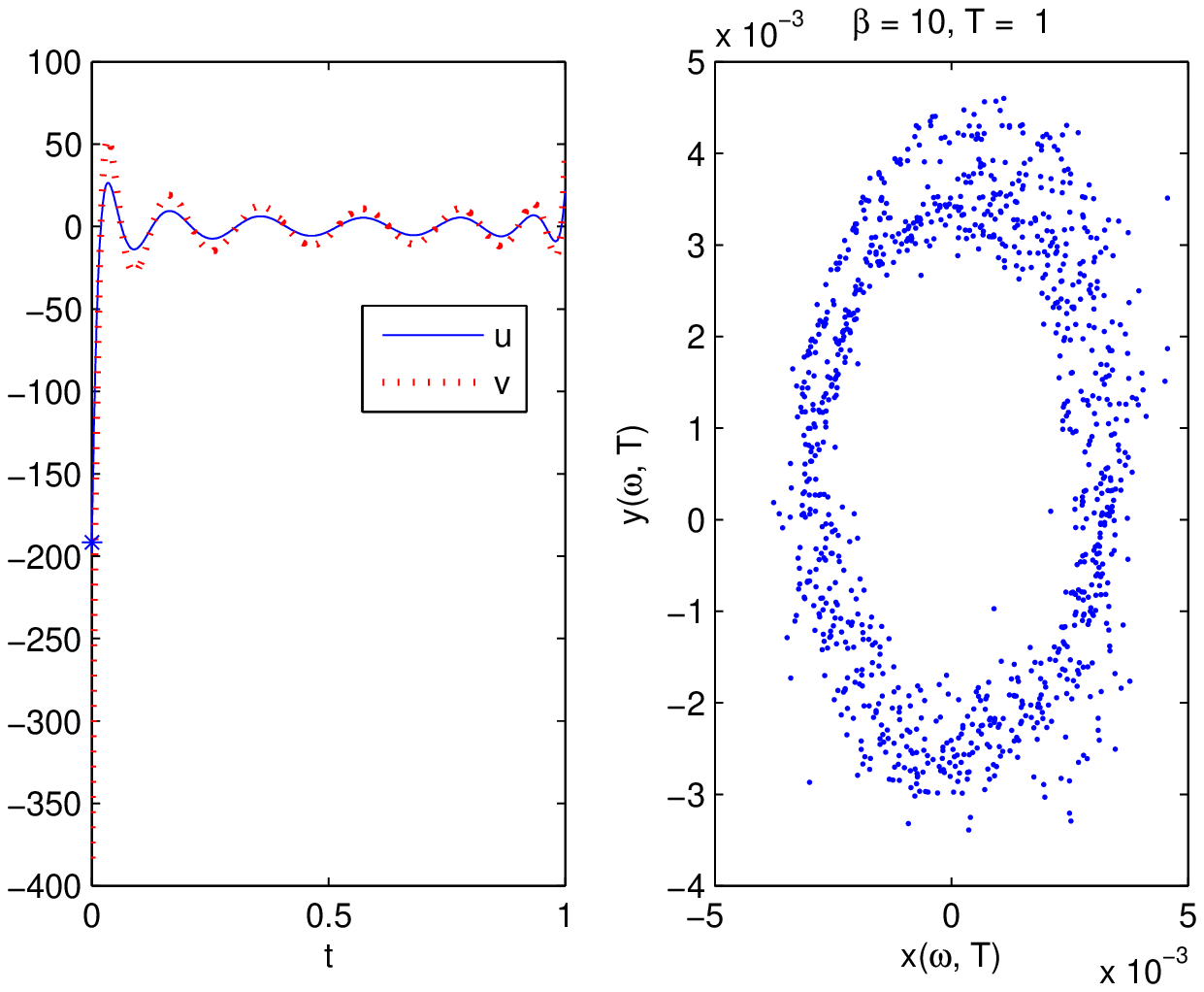}&\,
{\small (b)}\includegraphics[scale=0.5]{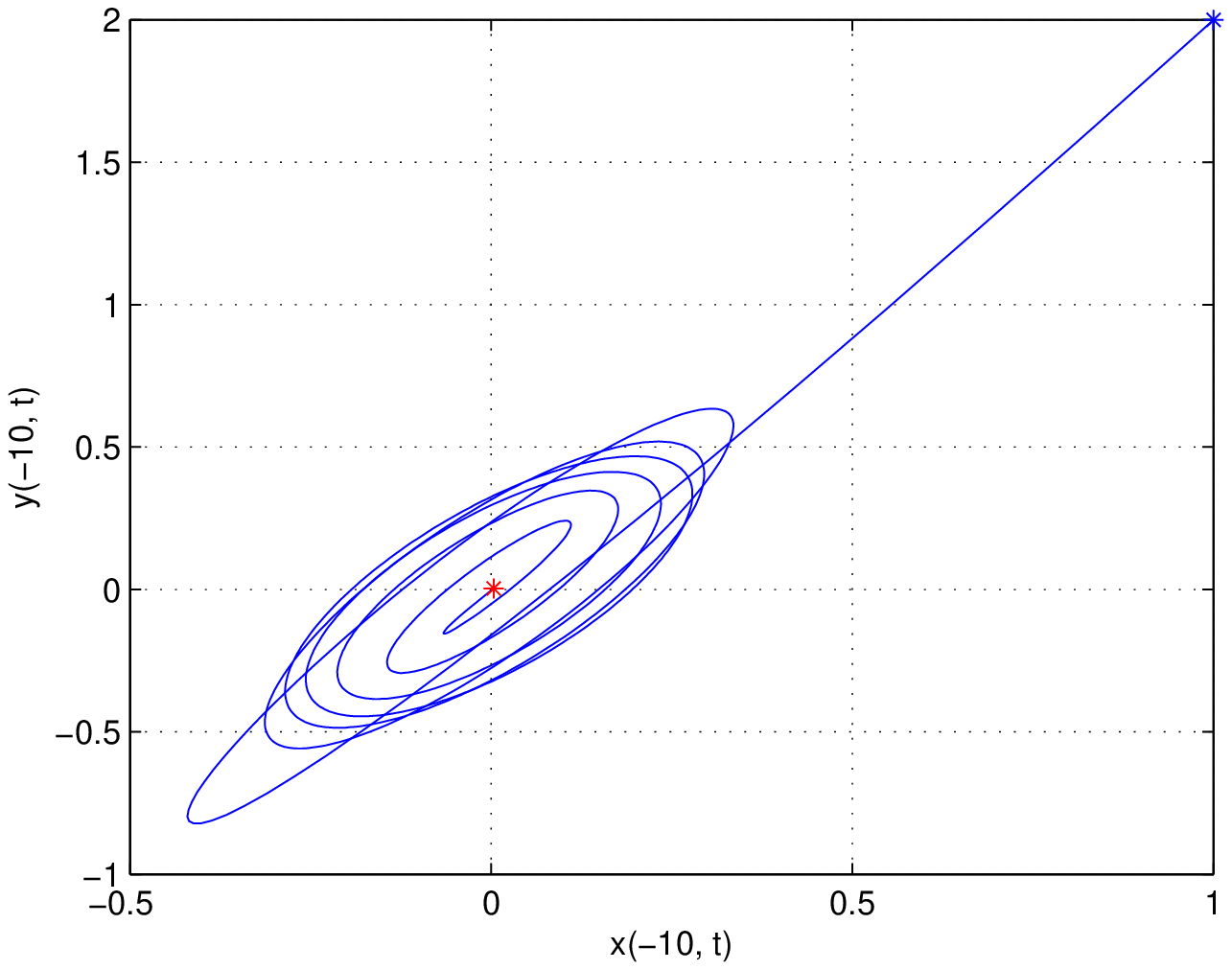}\\
{\small (c)}\includegraphics[scale=0.5]{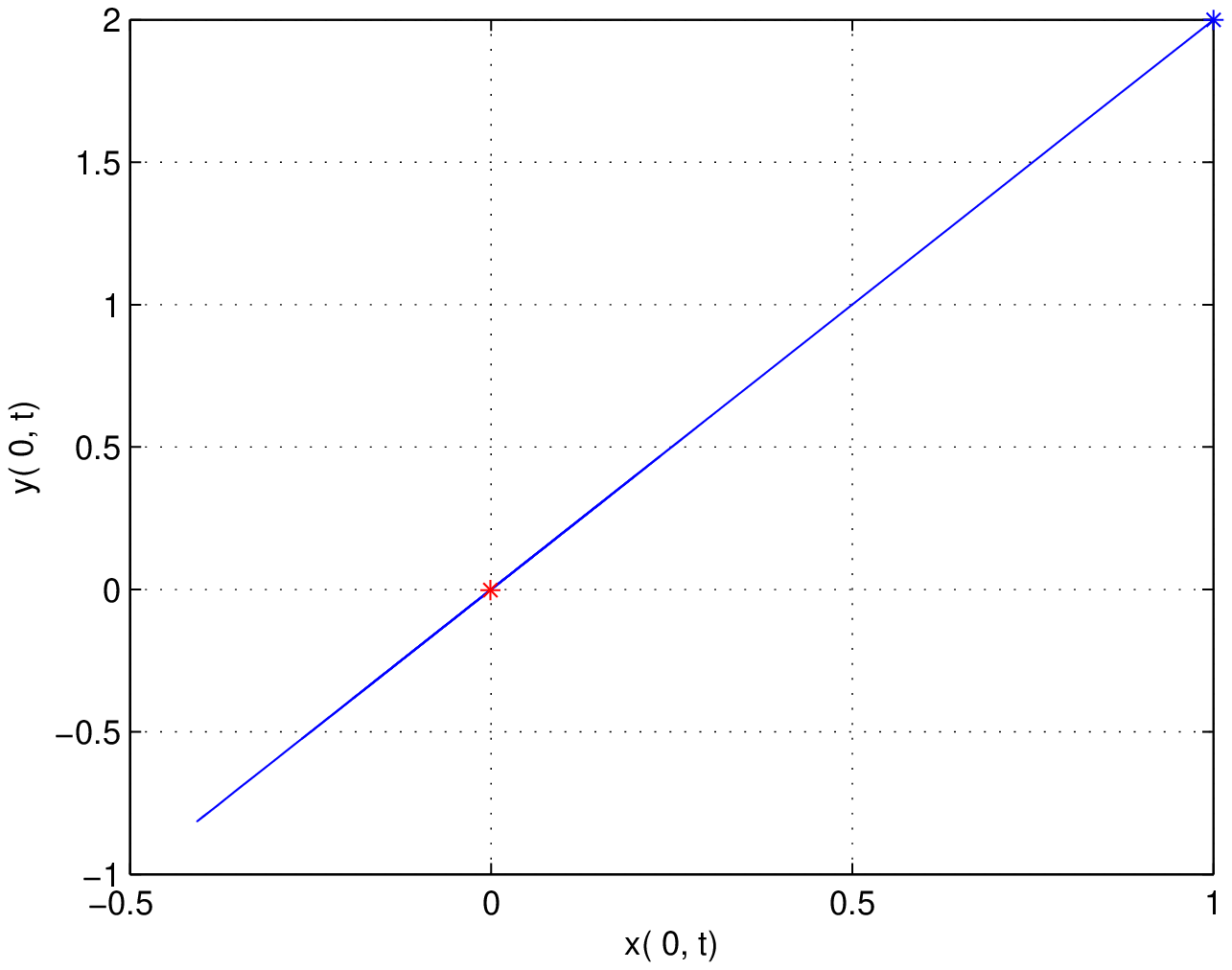}&\,
{\small (d)}\includegraphics[scale=0.5]{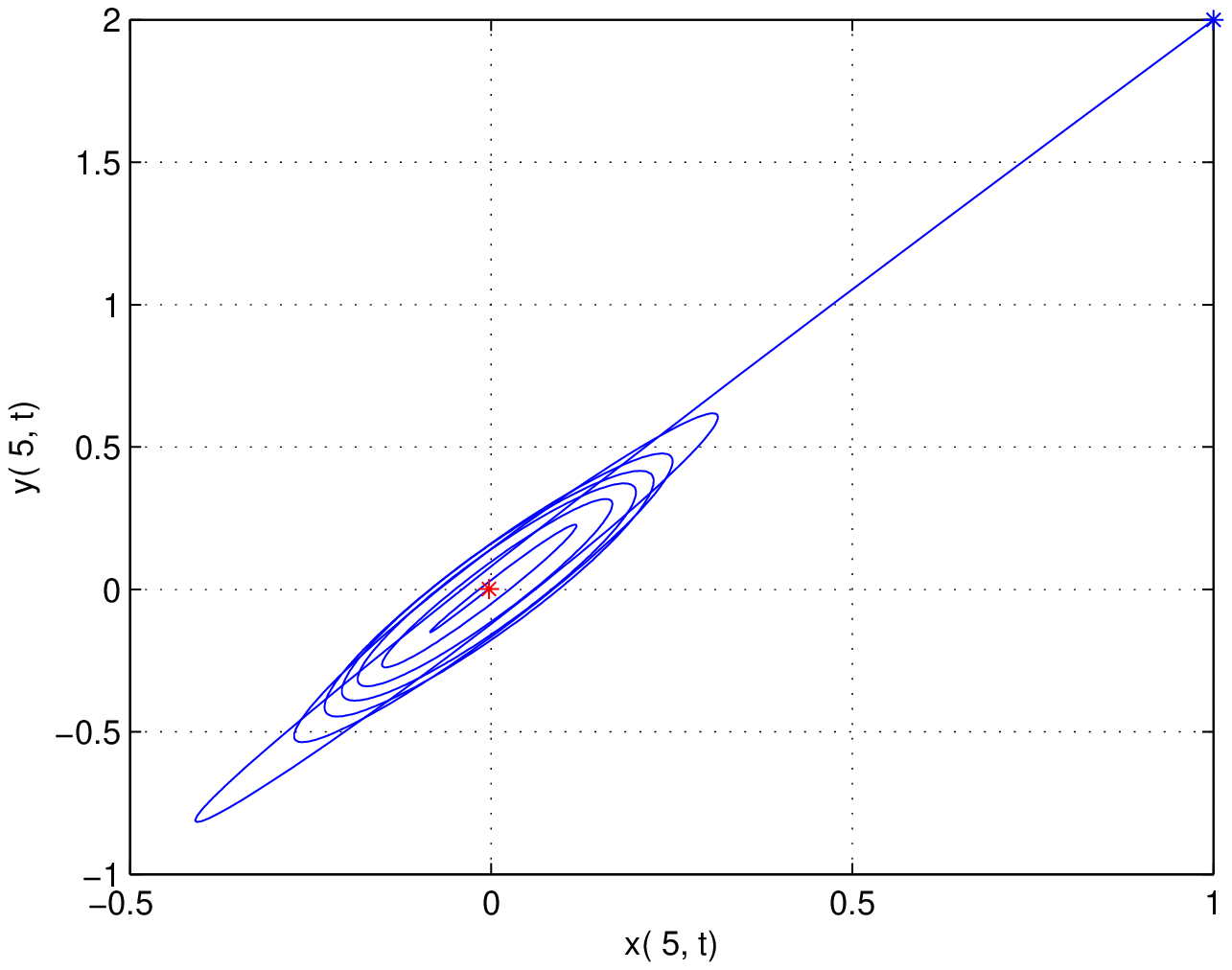}
\end{tabular}
\caption{\small The simulation results of Problem \ref{prob:harmonic1} for $N=1001$, $T=1$, and $\beta=10$. The initial state $X_0=(1,2)$ and the target state $X_F=(0,0)$. (a) The optimal control law $(u(t),v(t))$ for $t\in[0,1]$, and the final states for all systems $\omega\in[-10,10]$. (b) The trajectory for $\omega=-10$ following $(u(t),v(t))$. (c) The trajectory for $\omega=0$. (d) The trajectory for $\omega=5$.}\label{fig:case2}
\end{figure}

\subsection{Simulations}
\label{sec:simulation}
\begin{figure}[t]
\centering
\begin{tabular}{cc}
{\small (a)}\includegraphics[scale=0.52]{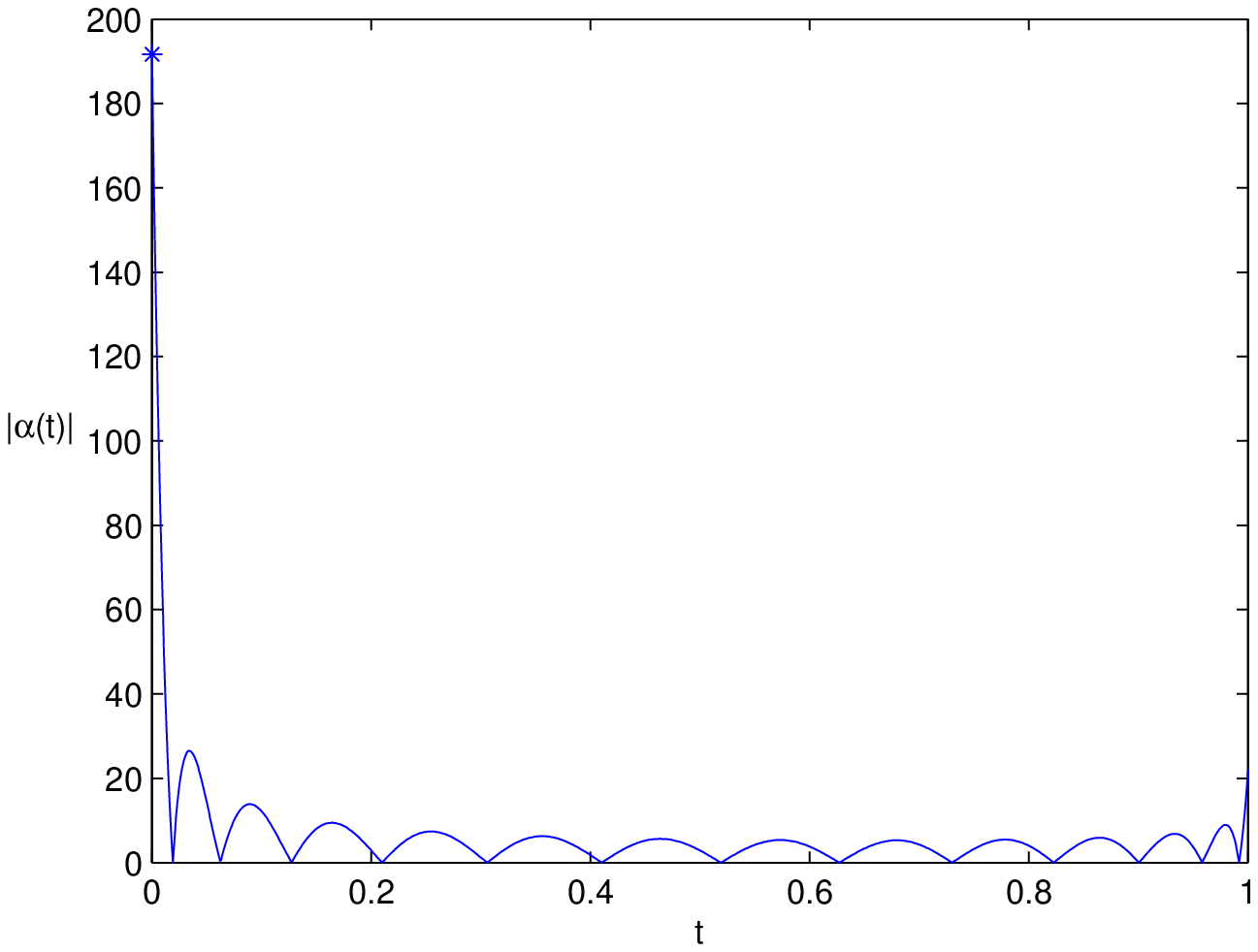}&\,
{\small (b)}\includegraphics[scale=0.52]{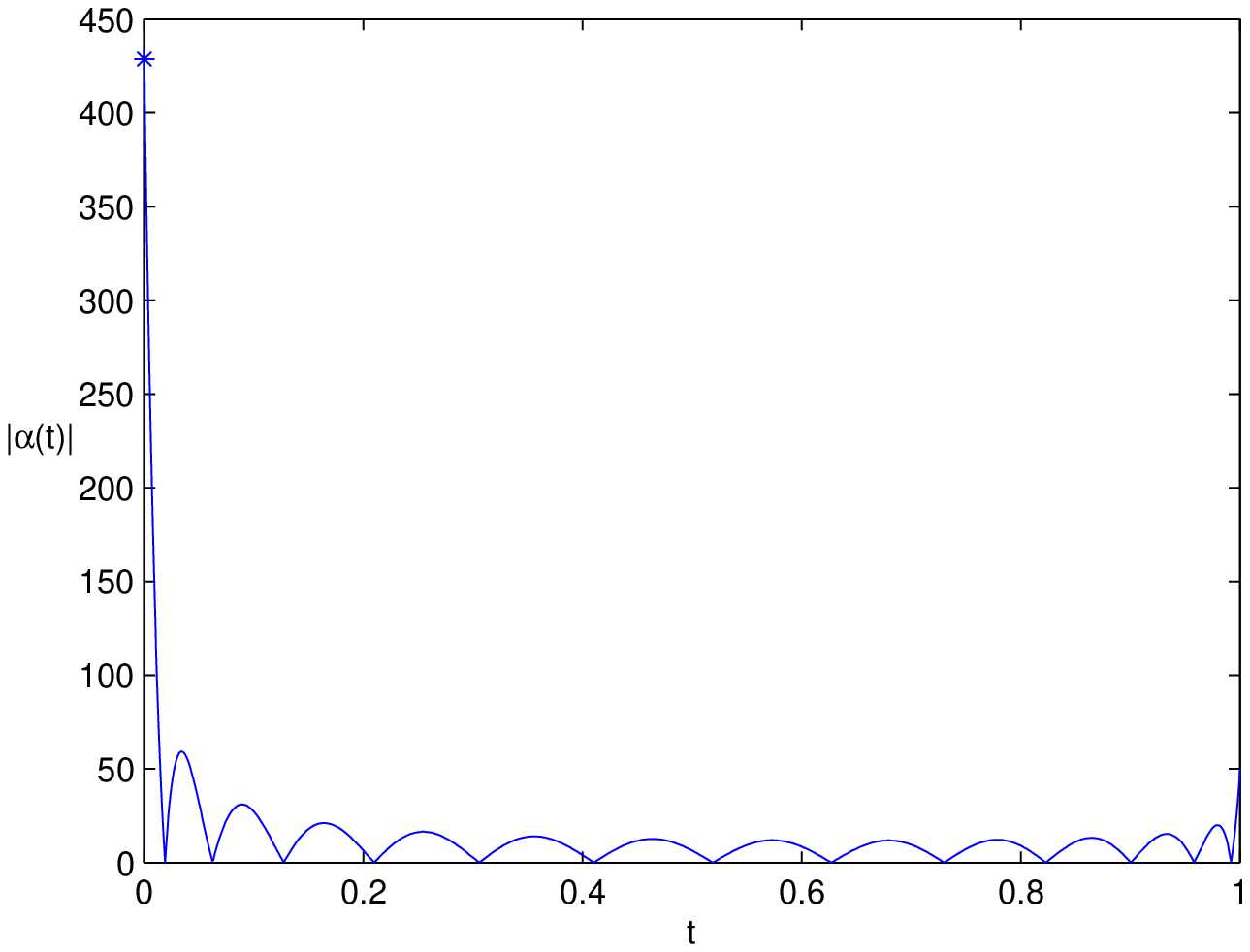}
\end{tabular}
\caption{\small (a) The control amplitude for the case (1). (b) The control amplitude for the case (2).}
\label{fig:amplitude}
\end{figure}

Here, we provide numerical solutions for $\alpha_N(t)$ since it is not of closed form. As shown in (\ref{eq:alphaN}), the ensemble control law $\alpha_N(t)$ is synthesized by the set of eigenfunctions $\{\tilde{\phi}_{n}\}$ and the corresponding eigenvalues $\{\lambda_{n}\}$ associated with the pswf's. These functions can be approximated by the \emph{discrete prolate spheroidal sequences (dpss's)}, denoted as $\{v_{t, k}(N, W)\}$, which are defined via the solution to the following equation \cite{Percival, Slepian4}
\begin{eqnarray}
\sum_{t^{\prime}=0}^{N-1}\frac{\sin\big[2\pi W(t-t^{\prime})\big]}{\pi(t-t^{\prime})}\, v_{t^{\prime}, k}(N,
W)=\lambda_{k}(N, W)v_{t, k}(N, W),\nonumber
\end{eqnarray}
where $0<W<\frac{1}{2}$ and $t = 0, 1, \ldots, N-1$. It is equivalent to saying that $\lambda_{k}(N, W)$ are the eigenvalues of the $N\times N$ matrix $A$ whose $(t, t^{\prime})$th element is
\begin{eqnarray}
\label{eq:dpssmatrix}(A)_{t,t^{\prime}}=\frac{\sin\big[2\pi W(t-t^{\prime})\big]}{\pi(t-t^{\prime})},\quad t,
t^{\prime}=0,1,\ldots,N-1,
\end{eqnarray}
and that the $N$ elements of the corresponding eigenvectors for this matrix are in fact subsequences of length $N$ of the \emph{dpss}'s. Note that $\lambda_{k}(N, W)$ are distinct, real, and ordered non-zero eigenvalues such that
\begin{eqnarray}
1>\lambda_{0}(N, W)>\lambda_{1}(N, W)>\ldots >\lambda_{N-1}(N,
W)>0,\nonumber
\end{eqnarray}
and the $dpss$'s are real-valued. Now, we show how to compute $\alpha_{N}(t)$. We present two cases with different initial states for $\beta=10$ and $T=1$:
\begin{enumerate}
\item[(1)] Consider $X(0,\omega)=(1,0)$ and $X(1,\omega)=(0,0)$. Then we have $p(0,\omega)=1$, $p(1,\omega)=0$, and hence, by (\ref{eq:varnat}), $\xi(\omega)=-1$ is a constant function.
\item[(2)] Consider $X(0,\omega)=(1,2)$ and $X(1,\omega)=(0,0)$. Then we have $p(0,\omega)=1+2i$, $p(1,\omega)=0$, and hence $\xi(\omega)=-1-2i$.
\end{enumerate}
According to the analysis above, the ``sinc'' kernel in (\ref{eq:pswf}) is replaced by the symmetric matrix $A$ as in (\ref{eq:dpssmatrix}), where $W=\frac{T\beta}{2\pi(N-1)}$. Note that the number of harmonic oscillators $N$ must be large enough to satisfy $W<\frac{1}{2}$ \cite{Percival}. Here we consider $N = 1001$ and the frequencies are uniformly sampled within $[-10,10]$. The simulation results are shown in Figure \ref{fig:case1} and Figure \ref{fig:case2}. It can be seen that following the resulting optimal control laws, the final states of all systems converge to a neighborhood of the desired target state, the origin. The trajectories for $\omega=-10$, $\omega=0$, and $\omega=5$ are displayed. Observe that in both cases a strong impulse is implemented initially, $|\alpha_N(0)|=191.7$ and $428.6$, respectively, as shown in Figure \ref{fig:amplitude}. Practical applications make it desirable to design a control with a limited amplitude. This leads to the following problem.

\subsection{Constrained Convex Optimization Problem}
In practice, the problem of interest is when the control amplitude is limited.

\begin{problem}
\label{prob:bounded}\rm Given a fixed time $T$, find bounded controls $u(t)$ and $v(t)$ satisfying the constraint $\sqrt{u^{2}(t)+v^{2}(t)}\leq A_{\rm{max}}\in\mathbb{R}^{+}$ for all $t\in[0,T]$, which will steer an ensemble of systems in (\ref{eq:harmonic}) from an initial state $X(0,\omega)=(1,0)$ as close as possible, in the $L_{2}$ sense, to the origin, $X(T,\omega)=(0,0)$, at time $T$.
\end{problem}

This problem can be formulated as the following minimization problem
\begin{eqnarray}
\label{eq:obj}
&& \min_{\alpha}\quad\int_{-\beta}^{\beta}\parallel p(T,\omega)-\mathbf{0}\parallel^2 d\omega, \\
\label{eq:cons} && {\rm s.t.}\quad\,\,\, u^2(t)+v^2(t)\leq A_{\rm{max}}^2,
\end{eqnarray}
where $p(T,\omega)$ defined in (\ref{eq:pomega}) depends on the control $\alpha$. It follows from (\ref{eq:pomega}) and (\ref{eq:obj}) that the problem can be further simplified as to minimize the following cost functional
\begin{eqnarray}
J=\int_{-\beta}^{\beta}\|\int_{0}^{T}e^{-i\omega\tau}\alpha(\tau)d\tau+1\|^2 d\omega,\nonumber
\end{eqnarray}
subject to the constraint (\ref{eq:cons}). By first integrating over $\omega$, we get
\begin{align}
J&=\int_{0}^{T}\int_{0}^{T}\frac{2\sin[\beta(\tau-\sigma)]}{\tau-\sigma}\,\alpha(\tau)\alpha^{\dag}(\sigma)d\tau
d\sigma\nonumber\\
&+\int_{0}^{T}\frac{2\sin(\beta\tau)}{\tau}\,\alpha(\tau)d\tau+\int_{0}^{T}\frac{2\sin(\beta\sigma)}{\sigma} \alpha^{\dag}(\sigma)d\sigma+2\beta.\nonumber
\end{align}
Observe that the imaginary part of the double integration in the expression for $J$ vanishes by antisymmetry, i.e.
\begin{eqnarray}
\int_{0}^{T}\int_{0}^{T}\frac{2\sin[\beta(\tau-\sigma)]}{\tau-\sigma} \Big[v(\tau)u(\sigma)-u(\tau)v(\sigma)\Big]d\tau d\sigma=0.\nonumber
\end{eqnarray}
Moreover, the sinc kernel is positive definite since
\begin{eqnarray}
\label{eq:PD}
\int_{0}^{T}\int_{0}^{T}\frac{2\sin[\beta(\tau-\sigma)]}{\tau-\sigma}\,
v(\tau)v(\sigma)d\tau\,d\sigma\nonumber\\
=\int_{-\beta}^{\beta}\parallel\int_{0}^{T}e^{-i\omega\tau}v(\tau)d\tau\parallel^2\,
d\omega > 0.
\end{eqnarray}
According to these observations, we can always minimize $J$ by the appropriate choice of $u(t)$ disregarding $v(t)$, because
\begin{align}
J & \geq\int_{0}^{T}\int_{0}^{T}\frac{2\sin[\beta(\tau-\sigma)]}{\tau-\sigma}u(\tau)u(\sigma)d\tau d\sigma\nonumber\\
& +\int_{0}^{T}\frac{4\sin(\beta\tau)}{\tau}u(\tau)d\tau+2\beta\nonumber.
\end{align}
Without loss of generality, we assume $A_{\rm{max}}=1$ and $\beta=1$ such that the ratio $\frac{A_{\max}}{\beta}=1$. The original problem described in (\ref{eq:obj}) and (\ref{eq:cons}) can now be recapitulated as follows:
\begin{align}
&\min_{u}\,\,\int_{0}^{T}\left[\int_{0}^{T}\frac{\sin(\tau-\sigma)}{\tau-\sigma}u(\tau)u(\sigma) d\sigma+\frac{2\sin(\tau)}{\tau}u(\tau)\right]d\tau\nonumber\\
\label{eq:consnew} & \ {\rm s.t.}\quad\,\, u^2(t)\leq 1.
\end{align}

\begin{prop}\label{prop:convex}
Let $\mathcal{S}=:\Big\{(u,v)\,\Big\vert\, u^2(t)+v^2(t)\leq A_{\rm{max}}^2\Big\}$. A local minimum of the cost functional $J$ over $\mathcal{S}$ is the global minimum.
\end{prop}
{\it Proof.} First, observe that $\mathcal{S}$ is a convex set. Furthermore, the cost function $J$ is quadratic in $\alpha$ with positive definite Hessian (see (\ref{eq:PD})). \hfill $\Box$

\begin{figure}[t]
\centering
\begin{tabular}{cc}
{\small (a)}\includegraphics[scale=0.44]{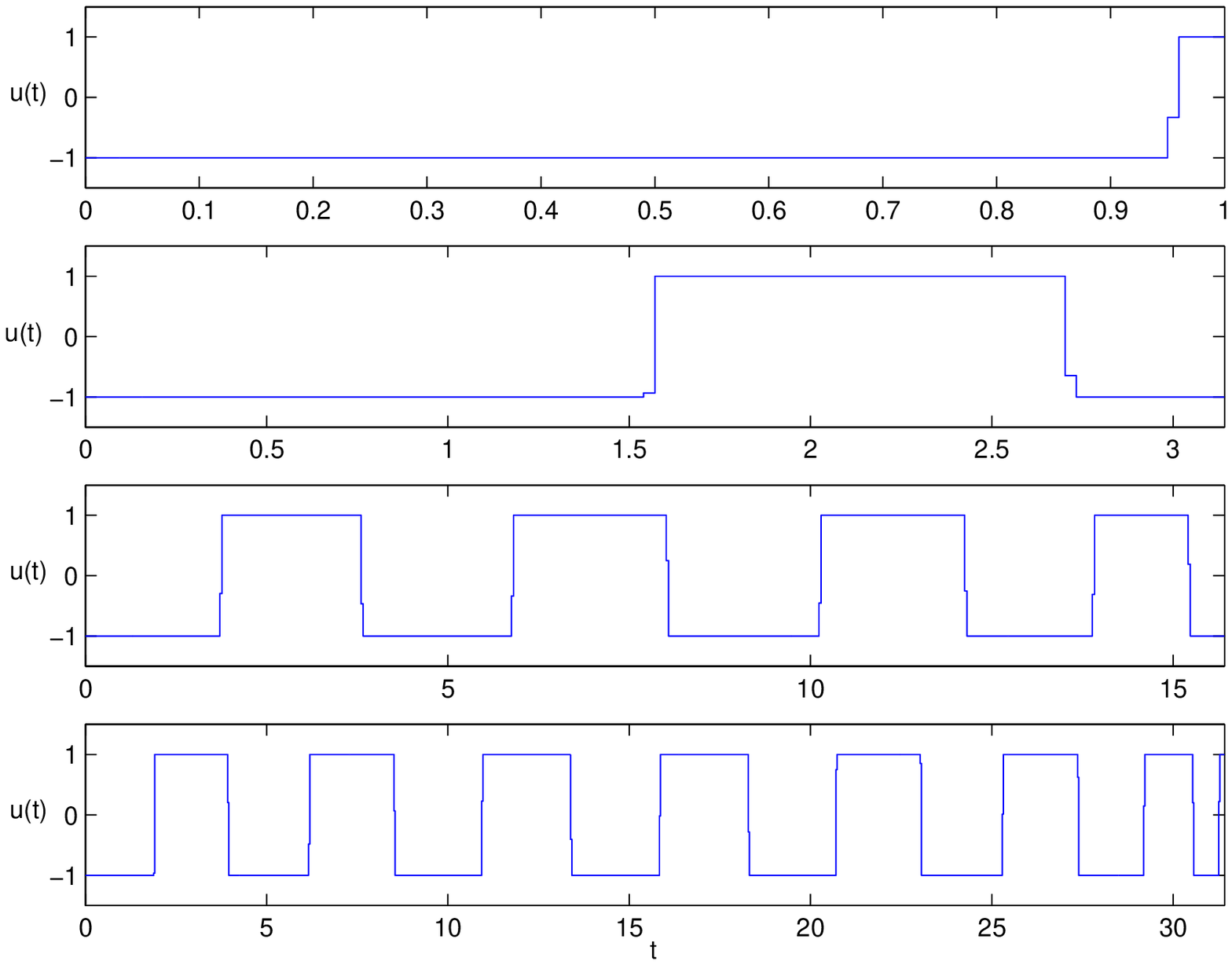}&
{\small (b)}\includegraphics[scale=0.44]{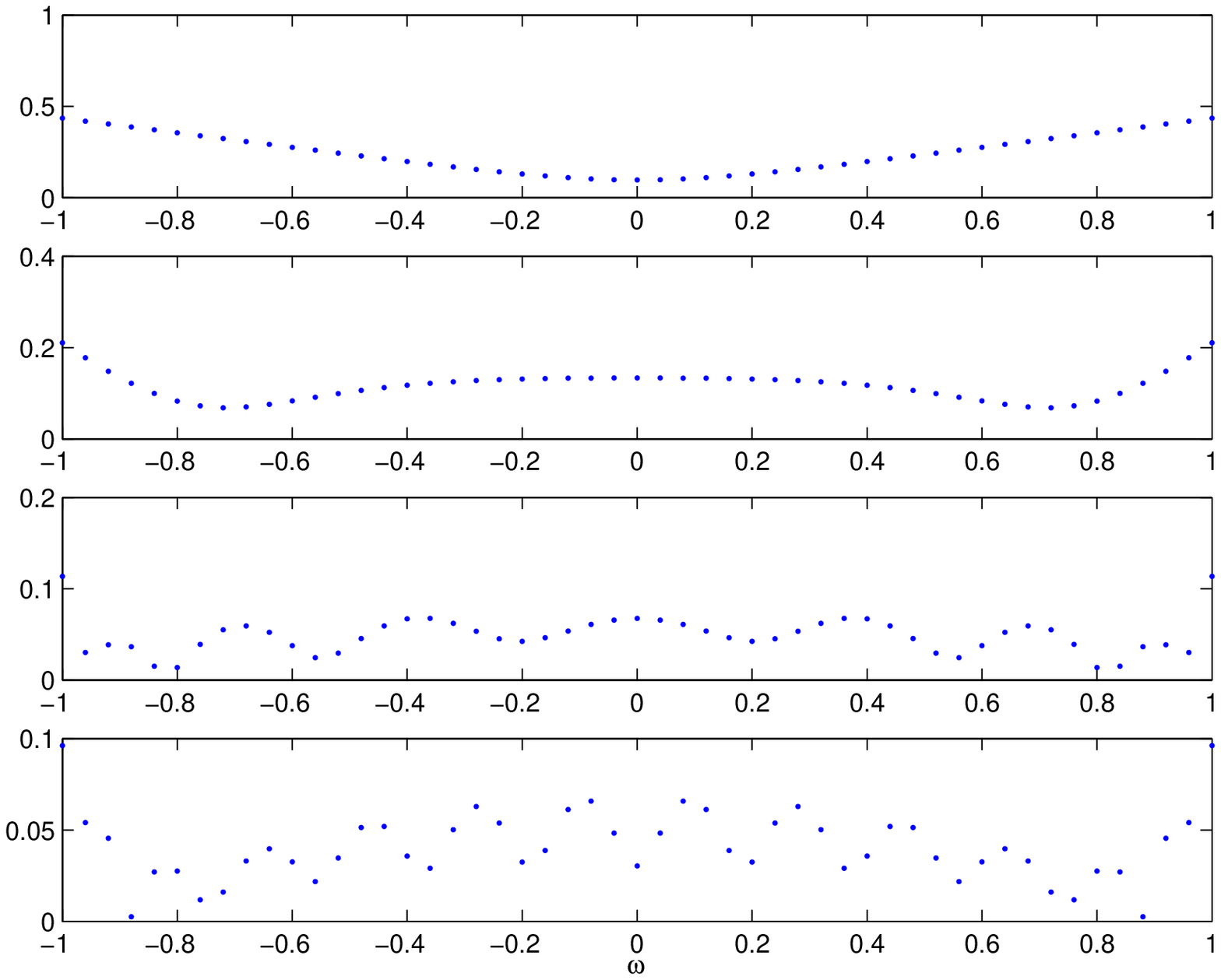}
\end{tabular}
\caption[The simulation results using Iterative Algorithm.]{(a) The optimal control laws for $T=1$, $T=\pi$, $T=5\pi$, and $T=10\pi$. (b) The distance between the final state and the origin, $|p(T,\omega)-{\bf 0}|$, of 51 harmonic oscillators for $T=1$, $T=\pi$, $T=5\pi$, and $T=10\pi$.}\label{fig:constrained}
\end{figure}

%

As a result, the model described in (\ref{eq:consnew}) is a convex optimization problem with a unique global minimum. We solve this problem numerically as a discrete quadratic optimization problem of the form
\begin{eqnarray}
\label{eq:objdis}
&&\min_{X}\quad X^{t}HX+2X^{t}Q \\
\label{eq:consdis} && {\rm s.t.}\quad\, |x_{i}|\leq 1, \quad i=1, 2,
\ldots, n,
\end{eqnarray}
where $X=(x_{1},x_{2}\ldots,x_{n})^{T}$, $t_{1}=0$, $t_{n}=T$, and
\begin{eqnarray}
H=\left[\begin{array}{cccc}\frac{\sin(t_{1}-t_{1})}{t_{1}-t_{1}} & \frac{\sin(t_{1}-t_{2})}{t_{1}-t_{2}} & \ldots & \frac{\sin(t_{1}-t_{n})}{t_{1}-t_{n}}\\ \frac{\sin(t_{2}-t_{1})}{t_{2}-t_{1}} & \frac{\sin(t_{2}-t_{2})}{t_{2}-t_{2}} & \ldots&\frac{\sin(t_{2}-t_{n})}{t_{2}-t_{n}}
\\ \vdots & \vdots & \ddots & \vdots\\ \frac{\sin(t_{n}-t_{1})}{t_{n}-t_{1}} & \frac{\sin(t_{n}-t_{2})}{t_{n}-t_{2}} & \ldots & \frac{\sin(t_{n}-t_{n})}{t_{n}-t_{n}}\end{array}\right],\,\,
Q=\left[\begin{array}{c}\frac{\sin(t_{1})}{t_{1}}\\ \frac{\sin(t_{2})}{t_{2}}\\ \vdots\\
\frac{\sin(t_{n})}{t_{n}}\end{array}\right], \nonumber
\end{eqnarray}
Some simulation results of Problem \ref{prob:bounded} for various values of $T$ are shown in Figure \ref{fig:constrained}, where we assume $A_{\rm{max}}=1$ and $\beta=1$. We consider 51 harmonic oscillators ($n=51$) with their frequencies uniformly distributed in $[-1,1]$. The optimal control laws with square wave forms are illustrated in Figure \ref{fig:constrained}(a). Figure \ref{fig:constrained}(b) shows the $L_2$ distance between the final state and the origin of each harmonic oscillator following the corresponding designed control laws.

\section{Conclusion}
In this paper, we studied ensemble control of general time-varying linear systems and derived the necessary and sufficient controllability conditions. The key idea of understanding controllability relies on investigating the solvability of the integral equation associated with the system dynamics, and it is of Fredholm equations of the first kind. We highlighted the role of singular systems and spectral theorem in designing ensemble control laws, and an analytical optimal control is provided. The work on computing singular values and eigenfunctions of operators will be pursued. Another interesting topic is to characterize the researchable set for the system \eqref{eq:linearTV} under the constrained controls, i.e., $u\in U\subset\mathbb{R}^{\rm m}$. We plan to continue working on extending and generalizing our current results towards the goal of developing a theory of ensemble control. We believe the study of ensemble control problems will foster further developments in control and systems theory with broader applications such as systems with parameter uncertainties as well as system identification.

\appendices
\section{}
\label{apd:Lcompact} We prove Proposition 1.

{\it Proposition 1:} The operator $L:\H_1\rightarrow\H_2$ defined by $$(Lu)(s)=\int_0^T\Phi(0,\tau;s)B(\tau,s)u(\tau)d\tau,$$ is compact, where $\Phi(t,0;s)$ satisfies for all $s\in[s_1,s_2]\subset\mathbb{R}$
$$\frac{d}{dt}\Phi(t,0;s)=A(t,s)\Phi(t,0;s);\quad\Phi(0,0;s)=I.$$
{\it Proof.} To prove this, we need the following tools.

\begin{defn}
Let $\H_0=L_2^{\rm n\times m}([s_1,s_2]\times [0,T])$ be the vector space of all those matrix valued functions $f$ whose elements $f_{ij}(s,t)$, $i=1,\ldots,n$, $j=1,\ldots,m$, are Lebesgue measurable on $[s_1,s_2]\times [0,T]$ and for which
$$\int_{s_1}^{s_2}\int_0^T\|f(s,t)\|^2 dt\,ds<\infty.$$
With the inner product for $f,g\in\H_0$ defined by
$$\< f,g\> =tr\int_{s_1}^{s_2}\int_0^T f(s,t)g(s,t)^{\dagger}\,dt ds,$$
$\H_0$ is a Hilbert space.
\end{defn}

Let $h(s,t)=\Phi(0,t;s)B(t,s)$, and then we have $$(Lu)(s)=\int_0^T h(s,t)u(t)dt.$$
We first show that $L$ is a bounded operator.

\begin{lem}
If $h(s,t)\in\H_0$, then $L\in\mathcal{B}(\H_1,\H_2)$.
\end{lem}
{\it Proof.} By Schwarz's inequality,
\begin{align*}
\|Lu\|^2 &= \int_{s_1}^{s_2}\left(\int_0^T h(s,t)u(t)dt\right)^{\dagger}\left(\int_0^T h(s,t)u(t)dt\right)ds\\
& \leq\|u\|^2\int_{s_1}^{s_2}\int_0^T h(s,t)^{\dagger}h(s,t)dtds\\
& \leq\|u\|^2\|h\|^2.
\end{align*}
Thus,
\begin{equation}
\label{eq:Lbdd}
\|L\|\leq\|h\|<\infty.\qquad\qquad\qquad\qquad\qquad\hfill\Box
\end{equation}

\begin{lem}
\label{lem:basis}
Suppose that $\{\phi_1,\phi_2,\ldots\}$ and $\{\psi_1,\psi_2,\ldots\}$ are orthonormal bases for $\H_1$ and $\H_2$, respectively, and
$$\Psi_{ij}(s,t)=\psi_i(s)\phi_j(t)^{\dagger}$$
for $(s,t)\in[s_1,s_2]\times [0,T], (i,j)\in\mathbb{N}\times\mathbb{N}$. Then, $\{\Psi_{ij}\}$ is an orthonormal basis of $\H_0$.
\end{lem}
{\it Proof.} For $j,k,m,n\in\mathbb{N}$,
\begin{align*}
\< \Psi_{jk},\Psi_{mn}\> &= tr\int_{s_1}^{s_2}\int_0^T\left[\psi_j(s)\phi_k(t)^{\dagger}\right] \left[\psi_m(s)\phi_n(t)^{\dagger}\right]^{\dagger} dtds\\
&= tr\int_{s_1}^{s_2}\psi_j(s)\left[\int_0^T\phi_k(t)^{\dagger}\phi_n(t)dt\right]\psi_m(s)^{\dagger}ds\\
&= \delta_{kn}\left[tr\int_{s_1}^{s_2}\psi_j(s)\psi_m(s)^{\dagger}ds\right]\\
&= \delta_{kn}\left[\int_{s_1}^{s_2}\psi_m(s)^{\dagger}\psi_j(s)ds\right]\\
&= \delta_{kn}\delta_{mj},
\end{align*}
so that $\{\Psi_{ij}\}$ is an orthonormal set in $\H_0$. 

Now, suppose that $k\in\H_0$ be such that $\< k,\Psi_{ij}\> =0$ for all $i,j\in\mathbb{N}$. Corresponding to this $k$, let $K:\H_1\rightarrow\H_2$ be defined by
$$(Ku)(s)=\int_0^T k(s,t)u(t)dt,\quad u\in\H_1.$$
Then, we have for all $i,j\in\mathbb{N}$
\begin{align*}
0=\< k,\Psi_{ij}\> &= tr\int_{s_1}^{s_2}\int_0^T k(s,t)\left[\psi_i(s)\phi_j(t)^{\dagger}\right]^{\dagger}dtds\\
&= tr\int_{s_1}^{s_2}\left[\int_0^T k(s,t)\phi_j(t)dt\right]\psi_i(s)^{\dagger}ds\\
&= tr\int_{s_1}^{s_2}(K\phi_j)\psi_i(s)^{\dagger}ds\\
&= \int_{s_1}^{s_2}\psi_i(s)^{\dagger}(K\phi_j)ds\\
&= \< \psi_i,K\phi_j\>_{\H_2}.
\end{align*}
Since $\psi_i$ is an orthonormal basis of $\H_2$ and $\phi_j$ is an orthonormal basis of $\H_1$, it follows that $$Ku=0 \quad\forall\, u\in\H_1.$$ Consequently, $k=0$ a.e. and hence $\{\Psi_{ij}\}$ is an orthonormal basis for $H_0$. $\hfill\Box$
%
\begin{thm}
\label{thm:compactseq}
Suppose $\{K_n\}$ is a sequence of compact operators in $\mathcal{B}(\H_1,\H_2)$ and $\|K_n-K\|\rightarrow 0$, where $K$ is in $\mathcal{B}(\H_1,\H_2)$. Then $K$ is compact.
\end{thm}
{\it Proof.} See \cite{Gohberg}. $\hfill\Box$

According to Lemma \ref{lem:basis}, $h(s,t)$ can be represented as $h=\sum_{i,j=1}^{\infty}\< h,\Psi_{ij}\> \Psi_{ij}$. Define
$$h_n(s,t)=\sum_{i,j=1}^{n}\< h,\Psi_{ij}\> \Psi_{ij}(s,t).$$
Then,
\begin{eqnarray}
\label{eq:h-hn}
\|h-h_n\|\rightarrow 0.
\end{eqnarray}
Let $L_n$ be the integral operator defined on $\H_1$ by
$$(L_n u)(s)=\int_0^T h_n(s,t)u(t)dt.$$
Now $L_n$ is a bounded linear operator of finite rank since $\mathcal{I}m L_n\subset\rm{sp}\{\psi_1,\ldots,\psi_n\}$, and hence $L_n$ is compact. By (\ref{eq:Lbdd}) and (\ref{eq:h-hn}) applied to $L-L_n$,
$$\|L-L_n\|\leq\|h-h_n\|\rightarrow 0.$$
This follows $L$ is compact by Theorem \ref{thm:compactseq}. Note that the compactness of $L$ can also be shown by using the properties of Hilbert-Schmidt operators \cite{Gohberg}. $\hfill\Box$

\end{document}